\expandafter\ifx\csname amssym.def\endcsname\relax \else\endinput\fi
%
\expandafter\edef\csname amssym.def\endcsname{%
       \catcode`\noexpand\@=\the\catcode`\@\space}
\catcode`\@=11
%

\def\undefine#1{\let#1\undefined}
\def\newsymbol#1#2#3#4#5{\let\next@\relax
 \ifnum#2=\@ne\let\next@\msafam@\else
 \ifnum#2=\tw@\let\next@\msbfam@\fi\fi
 \mathchardef#1="#3\next@#4#5}
\def\mathhexbox@#1#2#3{\relax
 \ifmmode\mathpalette{}{\m@th\mathchar"#1#2#3}%
 \else\leavevmode\hbox{$\m@th\mathchar"#1#2#3$}\fi}
\def\hexnumber@#1{\ifcase#1 0\or 1\or 2\or 3\or 4\or 5\or 6\or 7\or 8\or
 9\or A\or B\or C\or D\or E\or F\fi}

\font\tenmsa=msam10
\font\sevenmsa=msam7
\font\fivemsa=msam5
\newfam\msafam
\textfont\msafam=\tenmsa
\scriptfont\msafam=\sevenmsa
\scriptscriptfont\msafam=\fivemsa
\edef\msafam@{\hexnumber@\msafam}
\mathchardef\dabar@"0\msafam@39
\def\dashrightarrow{\mathrel{\dabar@\dabar@\mathchar"0\msafam@4B}}
\def\dashleftarrow{\mathrel{\mathchar"0\msafam@4C\dabar@\dabar@}}

\def\ulcorner{\delimiter"4\msafam@70\msafam@70 }
\def\urcorner{\delimiter"5\msafam@71\msafam@71 }
\def\llcorner{\delimiter"4\msafam@78\msafam@78 }
\def\lrcorner{\delimiter"5\msafam@79\msafam@79 }
\def\yen{{\mathhexbox@\msafam@55 }}
\def\checkmark{{\mathhexbox@\msafam@58 }}
\def\circledR{{\mathhexbox@\msafam@72 }}
\def\maltese{{\mathhexbox@\msafam@7A }}

\font\tenmsb=msbm10
\font\sevenmsb=msbm7
\font\fivemsb=msbm5
\newfam\msbfam
\textfont\msbfam=\tenmsb
\scriptfont\msbfam=\sevenmsb
\scriptscriptfont\msbfam=\fivemsb
\edef\msbfam@{\hexnumber@\msbfam}
\def\Bbb#1{{\fam\msbfam\relax#1}}
\def\widehat#1{\setbox\z@\hbox{$\m@th#1$}%
 \ifdim\wd\z@>\tw@ em\mathaccent"0\msbfam@5B{#1}%
 \else\mathaccent"0362{#1}\fi}
\def\widetilde#1{\setbox\z@\hbox{$\m@th#1$}%
 \ifdim\wd\z@>\tw@ em\mathaccent"0\msbfam@5D{#1}%
 \else\mathaccent"0365{#1}\fi}
\font\teneufm=eufm10
\font\seveneufm=eufm7
\font\fiveeufm=eufm5
\newfam\eufmfam
\textfont\eufmfam=\teneufm
\scriptfont\eufmfam=\seveneufm
\scriptscriptfont\eufmfam=\fiveeufm

\csname amssym.def\endcsname

\expandafter\ifx\csname pre amssym.tex at\endcsname\relax \else \endinput\fi
\expandafter\chardef\csname pre amssym.tex at\endcsname=\the\catcode`\@
\catcode`\@=11
\newsymbol\boxdot 1200
\newsymbol\boxplus 1201
\newsymbol\boxtimes 1202
\newsymbol\square 1003
\newsymbol\blacksquare 1004
\newsymbol\centerdot 1205
\newsymbol\lozenge 1006
\newsymbol\blacklozenge 1007
\newsymbol\circlearrowright 1308
\newsymbol\circlearrowleft 1309
\undefine\rightleftharpoons
\newsymbol\rightleftharpoons 130A
\newsymbol\leftrightharpoons 130B
\newsymbol\boxminus 120C
\newsymbol\Vdash 130D
\newsymbol\Vvdash 130E
\newsymbol\vDash 130F
\newsymbol\twoheadrightarrow 1310
\newsymbol\twoheadleftarrow 1311
\newsymbol\leftleftarrows 1312
\newsymbol\rightrightarrows 1313
\newsymbol\upuparrows 1314
\newsymbol\downdownarrows 1315
\newsymbol\upharpoonright 1316
 \let\restriction\upharpoonright
\newsymbol\downharpoonright 1317
\newsymbol\upharpoonleft 1318
\newsymbol\downharpoonleft 1319
\newsymbol\rightarrowtail 131A
\newsymbol\leftarrowtail 131B
\newsymbol\leftrightarrows 131C
\newsymbol\rightleftarrows 131D
\newsymbol\Lsh 131E
\newsymbol\Rsh 131F
\newsymbol\rightsquigarrow 1320
\newsymbol\leftrightsquigarrow 1321
\newsymbol\looparrowleft 1322
\newsymbol\looparrowright 1323
\newsymbol\circeq 1324
\newsymbol\succsim 1325
\newsymbol\gtrsim 1326
\newsymbol\gtrapprox 1327
\newsymbol\multimap 1328
\newsymbol\therefore 1329
\newsymbol\because 132A
\newsymbol\doteqdot 132B
 
\newsymbol\triangleq 132C
\newsymbol\precsim 132D
\newsymbol\lesssim 132E
\newsymbol\lessapprox 132F
\newsymbol\eqslantless 1330
\newsymbol\eqslantgtr 1331
\newsymbol\curlyeqprec 1332
\newsymbol\curlyeqsucc 1333
\newsymbol\preccurlyeq 1334
\newsymbol\leqq 1335
\newsymbol\leqslant 1336
\newsymbol\lessgtr 1337
\newsymbol\backprime 1038
\newsymbol\risingdotseq 133A
\newsymbol\fallingdotseq 133B
\newsymbol\succcurlyeq 133C
\newsymbol\geqq 133D
\newsymbol\geqslant 133E
\newsymbol\gtrless 133F
\newsymbol\sqsubset 1340
\newsymbol\sqsupset 1341
\newsymbol\vartriangleright 1342
\newsymbol\vartriangleleft 1343
\newsymbol\trianglerighteq 1344
\newsymbol\trianglelefteq 1345
\newsymbol\bigstar 1046
\newsymbol\between 1347
\newsymbol\blacktriangledown 1048
\newsymbol\blacktriangleright 1349
\newsymbol\blacktriangleleft 134A
\newsymbol\vartriangle 134D
\newsymbol\blacktriangle 104E
\newsymbol\triangledown 104F
\newsymbol\eqcirc 1350
\newsymbol\lesseqgtr 1351
\newsymbol\gtreqless 1352
\newsymbol\lesseqqgtr 1353
\newsymbol\gtreqqless 1354
\newsymbol\Rrightarrow 1356
\newsymbol\Lleftarrow 1357
\newsymbol\veebar 1259
\newsymbol\barwedge 125A
\newsymbol\doublebarwedge 125B
\undefine\angle
\newsymbol\angle 105C
\newsymbol\measuredangle 105D
\newsymbol\sphericalangle 105E
\newsymbol\varpropto 135F
\newsymbol\smallsmile 1360
\newsymbol\smallfrown 1361
\newsymbol\Subset 1362
\newsymbol\Supset 1363
\newsymbol\Cup 1264
 
\newsymbol\Cap 1265
 
\newsymbol\curlywedge 1266
\newsymbol\curlyvee 1267
\newsymbol\leftthreetimes 1268
\newsymbol\rightthreetimes 1269
\newsymbol\subseteqq 136A
\newsymbol\supseteqq 136B
\newsymbol\bumpeq 136C
\newsymbol\Bumpeq 136D
\newsymbol\lll 136E
 
\newsymbol\ggg 136F
 
\newsymbol\circledS 1073
\newsymbol\pitchfork 1374
\newsymbol\dotplus 1275
\newsymbol\backsim 1376
\newsymbol\backsimeq 1377
\newsymbol\complement 107B
\newsymbol\intercal 127C
\newsymbol\circledcirc 127D
\newsymbol\circledast 127E
\newsymbol\circleddash 127F
\newsymbol\lvertneqq 2300
\newsymbol\gvertneqq 2301
\newsymbol\nleq 2302
\newsymbol\ngeq 2303
\newsymbol\nless 2304
\newsymbol\ngtr 2305
\newsymbol\nprec 2306
\newsymbol\nsucc 2307
\newsymbol\lneqq 2308
\newsymbol\gneqq 2309
\newsymbol\nleqslant 230A
\newsymbol\ngeqslant 230B
\newsymbol\lneq 230C
\newsymbol\gneq 230D
\newsymbol\npreceq 230E
\newsymbol\nsucceq 230F
\newsymbol\precnsim 2310
\newsymbol\succnsim 2311
\newsymbol\lnsim 2312
\newsymbol\gnsim 2313
\newsymbol\nleqq 2314
\newsymbol\ngeqq 2315
\newsymbol\precneqq 2316
\newsymbol\succneqq 2317
\newsymbol\precnapprox 2318
\newsymbol\succnapprox 2319
\newsymbol\lnapprox 231A
\newsymbol\gnapprox 231B
\newsymbol\nsim 231C
\newsymbol\ncong 231D
\newsymbol\diagup 231E
\newsymbol\diagdown 231F
\newsymbol\varsubsetneq 2320
\newsymbol\varsupsetneq 2321
\newsymbol\nsubseteqq 2322
\newsymbol\nsupseteqq 2323
\newsymbol\subsetneqq 2324
\newsymbol\supsetneqq 2325
\newsymbol\varsubsetneqq 2326
\newsymbol\varsupsetneqq 2327
\newsymbol\subsetneq 2328
\newsymbol\supsetneq 2329
\newsymbol\nsubseteq 232A
\newsymbol\nsupseteq 232B
\newsymbol\nparallel 232C
\newsymbol\nmid 232D
\newsymbol\nshortmid 232E
\newsymbol\nshortparallel 232F
\newsymbol\nvdash 2330
\newsymbol\nVdash 2331
\newsymbol\nvDash 2332
\newsymbol\nVDash 2333
\newsymbol\ntrianglerighteq 2334
\newsymbol\ntrianglelefteq 2335
\newsymbol\ntriangleleft 2336
\newsymbol\ntriangleright 2337
\newsymbol\nleftarrow 2338
\newsymbol\nrightarrow 2339
\newsymbol\nLeftarrow 233A
\newsymbol\nRightarrow 233B
\newsymbol\nLeftrightarrow 233C
\newsymbol\nleftrightarrow 233D
\newsymbol\divideontimes 223E
\newsymbol\varnothing 203F
\newsymbol\nexists 2040
\newsymbol\Finv 2060
\newsymbol\Game 2061
\newsymbol\mho 2066
\newsymbol\eth 2067
\newsymbol\eqsim 2368
\newsymbol\beth 2069
\newsymbol\gimel 206A
\newsymbol\daleth 206B
\newsymbol\lessdot 236C
\newsymbol\gtrdot 236D
\newsymbol\ltimes 226E
\newsymbol\rtimes 226F
\newsymbol\shortmid 2370
\newsymbol\shortparallel 2371
\newsymbol\smallsetminus 2272
\newsymbol\thicksim 2373
\newsymbol\thickapprox 2374
\newsymbol\approxeq 2375
\newsymbol\succapprox 2376
\newsymbol\precapprox 2377
\newsymbol\curvearrowleft 2378
\newsymbol\curvearrowright 2379
\newsymbol\digamma 207A
\newsymbol\varkappa 207B
\newsymbol\Bbbk 207C
\newsymbol\hslash 207D
\undefine\hbar
\newsymbol\hbar 207E
\newsymbol\backepsilon 237F
\catcode`\@=\csname pre amssym.tex at\endcsname

%
%
\font\fivebi=cmmib5
\font\fivebsy=cmbsy5
\font\sixrm=cmr6
\font\sixi=cmmi6
\font\sixbf=cmbx6
\font\sixsy=cmsy6
\font\sixmsa=msam5 at 6pt
\font\sixmsb=msbm5 at 6pt
\font\sevenbi=cmmib7
\font\sevenbsy=cmbsy7
\font\eightrm=cmr8
\font\eightsl=cmsl8
\font\eightit=cmti8
\font\eighti=cmmi8
\font\eightbf=cmbx8
\font\eightsy=cmsy8
\font\eightmsa=msam7 at 8pt
\font\eightmsb=msbm7 at 8pt
\font\ninerm=cmr9
\font\ninesl=cmsl9
\font\nineit=cmti9
\font\ninei=cmmi9
\font\ninebf=cmbx9
\font\ninebi=cmmib10 scaled 900
\font\ninesy=cmsy9
\font\ninebsy=cmbsy10 scaled 900
\font\ninemsa=msam10 at 9pt
\font\ninemsb=msbm10 at 9pt
\font\tenbit=cmbxti10
\font\tenbsl=cmbxsl10
\font\tenbi=cmmib10
\font\tenbsy=cmbsy10
\font\twelvebf=cmbx12
\font\twelvebi=cmmib10 scaled 1200
\font\twelvebsy=cmbsy10 at 12pt

\let\sc=\sevenrm            
\def\eightpoint{%
     \def\rm{\fam0\eightrm}
     \textfont0=\eightrm \scriptfont0=\sixrm \scriptscriptfont0=\fiverm
     \textfont1=\eighti \scriptfont1=\sixi \scriptscriptfont1=\fivei
     \textfont2=\eightsy \scriptfont2=\sixsy \scriptscriptfont2=\fivesy
     \textfont3=\tenex \scriptfont3=\tenex \scriptscriptfont3=\tenex
     \textfont\itfam=\eightit \def\it{\fam\itfam\eightit}%
     \textfont\slfam=\eightsl \def\sl{\fam\slfam\eightsl}%
     \textfont\bffam=\eightbf \scriptfont\bffam=\sixbf
     \scriptscriptfont\bffam=\fivebf \def\bf{\fam\bffam\eightbf}%
     \textfont\msbfam=\eightmsb \textfont\msafam=\eightmsa
     \scriptfont\msafam=\sixmsa \scriptfont\msbfam=\sixmsb
     \scriptscriptfont\msafam=\fivemsa \scriptscriptfont\msbfam=\fivemsb
      \skewchar\eighti='177 \skewchar\sixi='177
      \skewchar\eightsy='60 \skewchar\sixsy='60
     \normalbaselineskip=10pt
     \setbox\strutbox=\hbox{\vrule height7pt depth3pt width0pt}
     \let\sc=\sixrm \let\big=\eightbig \normalbaselines\rm}
\def\ninepoint{%
     \def\rm{\fam0\ninerm}
     \textfont0=\ninerm \scriptfont0=\sixrm \scriptscriptfont0=\fiverm
     \textfont1=\ninei \scriptfont1=\sixi \scriptscriptfont1=\fivei
     \textfont2=\ninesy \scriptfont2=\sixsy \scriptscriptfont2=\fivesy
     \textfont3=\tenex \scriptfont3=\tenex \scriptscriptfont3=\tenex
     \textfont\itfam=\nineit \def\it{\fam\itfam\nineit}%
     \textfont\slfam=\ninesl \def\sl{\fam\slfam\ninesl}%
     \textfont\bffam=\ninebf \scriptfont\bffam=\sixbf
     \scriptscriptfont\bffam=\fivebf \def\bf{\fam\bffam\ninebf}%
     \textfont\msbfam=\ninemsb \textfont\msafam=\ninemsa
     \scriptfont\msafam=\sixmsa \scriptfont\msbfam=\sixmsb
     \scriptscriptfont\msafam=\fivemsa \scriptscriptfont\msbfam=\fivemsb
      \skewchar\ninei='177 \skewchar\sixi='177
      \skewchar\ninesy='60 \skewchar\sixsy='60
     \normalbaselineskip=11pt
     \setbox\strutbox=\hbox{\vrule height8pt depth3pt width0pt}%
     \let\sc=\sevenrm \let\big=\ninebig \normalbaselines\rm}
\def\tenpoint{%
     \def\rm{\fam0\tenrm}
     \textfont0=\tenrm \scriptfont0=\sevenrm \scriptscriptfont0=\fiverm
     \textfont1=\teni \scriptfont1=\seveni \scriptscriptfont1=\fivei
     \textfont2=\tensy \scriptfont2=\sevensy \scriptscriptfont2=\fivesy
     \textfont3=\tenex \scriptfont3=\tenex \scriptscriptfont3=\tenex
     \textfont\itfam=\tenit \def\it{\fam\itfam\tenit}%
     \textfont\slfam=\tensl \def\sl{\fam\slfam\tensl}%
     \textfont\bffam=\tenbf \scriptfont\bffam=\sevenbf
     \scriptscriptfont\bffam=\fivebf \def\bf{\fam\bffam\tenbf}%
     \textfont\msbfam=\tenmsb \textfont\msafam=\tenmsa
     \scriptfont\msafam=\sevenmsa \scriptfont\msbfam=\sevenmsb
     \scriptscriptfont\msafam=\fivemsa \scriptscriptfont\msbfam=\fivemsb
     \normalbaselineskip=12pt
     \let\sc=\sevenrm \let\big=\tenbig \normalbaselines\rm}
\catcode`@=11        
\def\tenbig#1{{\hbox{$\left#1\vbox to8.5pt{}\right.\n@space$}}}
\def\ninebig#1{{\hbox{$\textfont0=\tenrm\textfont2=\tensy
      \left#1\vbox to7.25pt{}\right.\n@space$}}}
\def\eightbig#1{{\hbox{$\textfont0=\ninerm\textfont2=\ninesy
      \left#1\vbox to6.5pt{}\right.\n@space$}}}
\catcode`@=12        
\def\bold{%
     \textfont0=\tenbf \scriptfont0=\sevenbf \scriptscriptfont0=\fivebf
     \textfont1=\tenbi \scriptfont1=\sevenbi \scriptscriptfont1=\fivebi
     \textfont2=\tenbsy \scriptfont2=\sevenbsy \scriptscriptfont2=\fivebsy
       \textfont\itfam=\tenbit \def\it{\fam\itfam\tenbit}%
       \textfont\slfam=\tenbsl \def\sl{\fam\slfam\tenbsl}%
       \textfont\bffam=\tenbf \scriptfont\bffam=\sevenbf
       \textfont\msbfam=\tenmsb \textfont\msafam=\tenmsa
   \fam0\tenbf}
\def\bigbold{%
     \textfont0=\twelvebf \scriptfont0=\ninebf \scriptscriptfont0=\sevenbf
     \textfont1=\twelvebi \scriptfont1=\ninebi \scriptscriptfont1=\sevenbi
     \textfont2=\twelvebsy \scriptfont2=\ninebsy \scriptscriptfont2=\sevenbsy
     \fam0\twelvebf}
%
%
%
%
\let\plainitem=\item
\let\plainitemitem=\itemitem
%

%
%
\def\openface{\Bbb}                
\def\N{{\openface N}}              

\def\C{{\openface C}}
%
%
%
\def\g{\hskip.17em\relax}               
\def\th{\thinspace}                     
\def\nl{\hfil\break}
\newskip\Bigskipamount
   \Bigskipamount=2\baselineskip plus.5\baselineskip minus.3\baselineskip
\def\Bigbreak{\removelastskip\vskip0pt plus .1\vsize\penalty-1000
              \vskip0pt plus-.1\vsize\vskip\Bigskipamount}
\def\Nobreak$$#1$${\postdisplaypenalty=10000$$#1$$\postdisplaypenalty=0}
%
%
%
\let\vvvvv=\v 
   
   \def\SC{Stone-\vvvvv Cech}
\let\doublebar=\| 
\def\|{\!\!\restriction\!\!}

  \let\sub=\sube

\def\supe{\supseteq}

\def\sm{\smallsetminus}
\def\es{\emptyset}

\def\T{T} 

\def\:{\colon}
\def\minor{\preccurlyeq} 
\def\Minor{\succcurlyeq}

\def\slt{\mathrel{\hbox{$\minor$\kern-.6em\lower.33ex\hbox{${}_s\;$}}}}
\def\sgt{\mathrel{\mathchoice                        
   {\hbox{$\Minor$\kern-.5em\lower.3ex\hbox{${}_s$}}}
   {\hbox{$\Minor$\kern-.5em\lower.3ex\hbox{${}_s$}}}
   {\hbox{$\scriptstyle\Minor\kern-.43em\lower.28ex\hbox{$\scriptstyle{}_s$}$}}
 {\hbox{$\scriptstyle\Minor\kern-.43em\lower.28ex\hbox{$\scriptstyle{}_s$}$}} }}    

\def\ucl(#1){\lfloor #1 \rfloor}
\def\dcl(#1){\lceil #1 \rceil}
\def\interior{\mathaccent"7017\relax}
%
%
\def\specrel#1#2{\mathrel{\mathop{\kern0pt #1}\limits_{#2}}}
\def\Specrel#1#2{\mathrel{\mathop{\kern0pt #1}\limits^{#2}}}
%
%
\def\alignspecrel#1#2{\mathrel{\mathop{\kern0pt #1}\limits_{\hbox
   to0pt{\hss$\scriptstyle#2$\hss}}}}
\def\alignSpecrel#1#2{\mathrel{\mathop{\kern0pt #1}\limits^{\hbox
   to0pt{\hss$\scriptstyle#2$\hss}}}}
\def\invlim{\specrel\lim{\raise 2pt\hbox{$\longleftarrow$}}}
\def\proof{\removelastskip\penalty55\medskip\noindent{\bf Proof. }}
\def\noproof{{\unskip\nobreak\hfill\penalty50\hskip2em\hbox{}\nobreak\hfill%
       $\square$\parfillskip=0pt\finalhyphendemerits=0\par}\goodbreak}
\def\endproof{\noproof\bigskip}

%
\newcount\refno
\def\ref#1#2\par{{\plainitem{[??]}#2\smallskip}}
\newtoks\thingtowrite 
\long\def\writerefnumber#1{%
    \thingtowrite={#1}%
    \immediate\write\refnumbersfile{\the\thingtowrite}%
    }
\newwrite\refnumbersfile
\def\makerefnumbers{\immediate\openout\refnumbersfile=RefNumbers%
  \refno=0 \writerefnumber{\refno=0}
  \def\ref##1##2\par{\global\advance\refno by 1
    \writerefnumber{\global\advance\refno by 1 \newcounter##1 ##1=\the\refno}%
    \plainitem{[\the\refno]}##2\smallskip
    }%
  }
\def\autorefnumbers{\refno=0
  \def\ref##1##2\par{\advance\refno by 1\plainitem{[\the\refno]}##2\smallskip}
  }
\def\userefnumbers{\refno=0
  \def\ref##1##2\par{\advance\refno by 1\plainitem{[\the##1]}##2\smallskip}
  }
%
%
\def\proclaimwithname #1. (#2) #3\par{{\bigbreak
  \clubpenalty=10000\noindent{\bf#1.\enspace}(#2)\nl
  {\sl #3}\par\bigbreak}}
\def\proposition (#1) #2\par{{\setbox0\hbox{(#1)\enspace}\bigbreak
   \sl\hangindent\the\wd0 \noindent\hskip\the\wd0
   \llap{\box0}\ignorespaces#2\par\bigbreak}}
\def\subsection #1\par{\vskip 3\medskipamount minus \smallskipamount\leftline{\bold #1}
        \penalty10000\smallskip\noindent}
      \def\section #1\par{\Bigbreak\centerline{\bf #1} 
              \penalty10000\bigskip\noindent}
%
\def\beginpsection #1\par{\Bigbreak\centerline{\bold #1}
        \penalty10000\bigskip}
\def\psubsection #1\par{\bigbreak\leftline{\bold #1}\penalty10000\bigskip}
%

%
%
\def\pitem#1{\smallskip\advance\parindent by 3mm
             \plainitem{\rm(#1)}\advance\parindent by-3mm}
\def\pitemitem#1{\smallskip\advance\parindent by 3mm
             \plainitemitem{\rm(#1)}\advance\parindent by-3mm}
\def\varitemitem#1{{\setbox0\hbox{\hskip\parindent#1\enskip}
           \smallbreak\hangindent\the\wd0 \noindent\hskip\the\wd0
           \llap{#1\enskip}\ignorespaces}}
%
%
\newdimen\newparindent
\def\iitem#1#2\par{\newparindent=\parindent \advance\newparindent by 3mm
           \smallbreak \hangindent\newparindent \noindent\hskip\newparindent
           \llap{{\rm #1}\enspace}\ignorespaces#2\par\smallbreak}
\def\iitemitem#1#2\par{\newparindent=\parindent \advance\newparindent by 3mm
           \smallbreak \hangindent2\newparindent \noindent\hskip2\newparindent
           \llap{{\rm #1}\enspace}\ignorespaces#2\par\smallbreak}
\def\varitem#1#2\par{{\setbox0\hbox{{\rm #1}\enspace}
           \smallbreak \hangindent\the\wd0 \noindent\hskip\the\wd0
           \llap{{\rm #1}\enspace}\ignorespaces#2\par\smallbreak}}
\def\enditem{\par}
%
\def\Textindent#1{\par \advance\parindent by 3mm
                  \textindent{{\rm #1}} \advance\parindent by -3mm}
\def\indentedline#1{\advance\hsize by -\parindent \line{#1}
                   \advance\hsize by \parindent}
\def\iindentedline#1{\advance\parindent by 3mm
                     \advance\hsize by -\parindent
                     \line{#1}
                     \advance\hsize by \parindent
                     \advance\parindent by -3mm}
\newdimen\margin   
\def\textdisplay#1&#2&#3$${\margin=\hsize
          \setbox1=\hbox{$\displaystyle#1\quad$}%
          \setbox2=\hbox{\quad#2\qquad$#3$}%
                     \advance\margin by-\wd1
                     \divide\margin by 2
   \ifdim\wd2 < \margin
      \box1\rlap{\quad#2}\eqno#3$$%
   \else
      \line{\qquad\hfil \box1\quad #2 \qquad $#3$}$$%
   \fi}
%
\def\ltextdisplay#1&#2&#3$${\margin=\hsize
           \setbox2=\hbox{$\displaystyle#2\quad$}
           \setbox3=\hbox{\quad#3\qquad}
                     \advance\margin by-\wd2
                     \divide\margin by 2
   \ifdim\wd3 < \margin
      \line{$#1$\hfil\box2\hbox to \margin{\box3\hfil}}$$%
   \else
      \line{$#1$\qquad\hfil\box2\quad #3\qquad} $$%
   \fi}
%
\def\textno#1&#2\par{%
   \margin=\hsize
   \advance\margin by -4\parindent
          \setbox1=\hbox{\sl#1}%
   \ifdim\wd1 < \margin
      $$\box1\eqno#2$$\endgraf%
   \else
      \bigbreak
      \line{\indent$\vcenter{\advance\hsize by -3\parindent
      \sl\noindent#1}\hfil#2$}%
      \bigbreak
   \fi}
%
\def\textlno#1&#2\par{%
   \margin=\hsize
   \advance\margin by -4\parindent
          \setbox1=\hbox{\sl#1}%
   \ifdim\wd1 < \margin
      $$\box1\leqno#2$$%
   \else
      \bigbreak
      \line{$#2\hfil\vcenter{\advance\hsize by -3\parindent
          \sl\noindent#1}\hskip\parindent$}%
      \bigbreak
   \fi}
%
%
%
\newcount\commentno
\def\COMMENT{}%

\def\nocomments{}
%
%
\def\?#1{\vadjust{\vbox to 0pt{\vss\vskip-8pt\leftline{%
     \llap{\hbox{\vbox{\pretolerance=-1
     \doublehyphendemerits=0\finalhyphendemerits=0
     \hsize16truemm\tolerance=10000\eightpoint
     \lineskip=0pt\lineskiplimit=0pt
     \rightskip=0pt plus16truemm\baselineskip8pt\noindent
     \hskip0pt        
     #1\endgraf}\hskip7truemm}}}\vss}}}
\def\d{}
%
%
%
%
\def\ds#1{}
%
%
\long\def\indexwrite#1{%
    \thingtowrite={#1}%
    \immediate\write\index{\the\thingtowrite}%
    }
%
%
\newwrite\index
\def\makeindex{\immediate\openout\index=index%
   \immediate\write\index{\catcode`@=11}%
   \def\d##1 {\ifmmode
     \write\index{$##1$, }%
     \write\index{\the\count0}\write\index{}
   \else
     \write\index{{##1}, }%
     \write\index{\the\count0}\write\index{}
   \fi {##1} }
      \def\ds##1{\ifmmode
     \write\index{$##1$, }%
     \write\index{\the\count0}\write\index{}
   \else
     \write\index{##1, }%
     \write\index{\the\count0}\write\index{}
   \fi}}
%
\newdimen\gap
\gap=3truemm
\newdimen\hackwidth
\hackwidth=15truemm
\def\disablems{\def\mos##1{\strut}}
\def\mo#1{\ifmmode {#1}\else {\it#1}\fi\mos{#1}}
\def\mos#1{\ifmmode
     \strut\vadjust{\vbox to 0pt{\vss\kern-11pt\leftline{%
     \llap{\hbox{\vbox{\pretolerance=-1
     \doublehyphendemerits=0\finalhyphendemerits=0
     \hsize\hackwidth\tolerance=10000\eightpoint
     \lineskip=0pt\lineskiplimit=0pt
     \rightskip=0pt plus\hsize\baselineskip8pt\noindent
     $#1$\strut\endgraf}\hskip\gap }}}\vss}}%
   \else
     \strut\vadjust{\vbox to 0pt{\vss\kern-11pt\leftline{%
     \llap{\hbox{\vbox{\pretolerance=-1
     \doublehyphendemerits=0\finalhyphendemerits=0
     \hsize\hackwidth\tolerance=10000\eightpoint
     \lineskip=0pt\lineskiplimit=0pt
     \rightskip=0pt plus\hsize\baselineskip8pt\noindent
     \hskip0pt    
     {\sl#1}\strut\endgraf}\hskip\gap }}}\vss}}%
   \fi}%
\newcount\remarkno
\def\REMARK#1{{\footnote{${}^{\the\remarkno}$}{{#1}}%
   \global\advance\remarkno by1}}
\def\noremarks{\def\REMARK##1{}}
%
%
%
%
\def\picture #1 by #2 (#3){
  \vbox to #2{
          \vfill
          \special{picture #3}
          \hrule width #1 height 0pt depth 0pt
           }}
\newdimen\topfiguremargin
   \topfiguremargin=0pt                                  
\newdimen\bottomfiguremargin
   \bottomfiguremargin=\medskipamount                    
\newdimen\normalpictureheight
\normalpictureheight=40mm
\def\Fig.#1 (#2by#3; heightfactor:#4; caption:#5) {{%
   \dimen2=\normalpictureheight
   \dimen0=#2                          
      \divide\dimen2 by 1000
      \multiply\dimen2 by#4              
   \count2=\dimen2                  
      \dimen1=#3                             
   \count1=\dimen1
   \divide\count1 by 1000
   \divide\count2 by \count1          
   \divide\dimen0 by 1000
   \multiply\dimen0 by \count2      
         \dimen1=\hsize
         \advance\dimen1 by -\dimen0
         \divide\dimen1 by 2               
   \midinsert
   \vbox to \topfiguremargin{\vfil}
   \noindent\hskip\dimen1
   \picture\dimen0 by \dimen2  (Fig.#1 scaled \the\count2)%
   \vskip\bottomfiguremargin                     
      \ninepoint
      \parindent=.1\hsize\narrower\narrower
      \setbox0\hbox{#5}
      \ifdim\wd0 < .6\hsize
           \centerline{F{\sc IGURE} #1.\hskip1em#5}
       \else
           \plainitem{F{\sc IGURE} #1. }#5\par
       \fi
   \vskip0pt\endinsert}}
%
\def\textpicture #1(#2by#3; #4width#5lower#6){{%
      \dimen0=#5\count2=\dimen0                    
      \dimen0=#2\count1=\dimen0                    
   \divide\count1 by 1000
   \divide\count2 by \count1                 
   \hbox{\vrule #4width0pt\vbox to 0pt{\vss\vskip#6%
      \special{picture #1 scaled \the\count2}\hrule width#5 height0pt\vss}}}}
%
%
%
%
%
\def\figure #1. #2 (#3; #4) {{%
   \def\bigskip{\par\ifdim\lastskip<\bigskipamount\removelastskip
      \vskip\bigskipamount\fi}
   \midinsert\vskip\topfiguremargin
   \dimen0=\normalpictureheight
      \divide\dimen0 by 1000
      \multiply\dimen0 by#4        
   \centerline{\epsfbox{#3.eps}}
   \vskip\bottomfiguremargin                     
      \ninepoint
      \parindent=.1\hsize\narrower\narrower
      \setbox0\hbox{#2}
      \ifdim\wd0 < .6\hsize
           \centerline{F{\sc IGURE} #1.\hskip1em#2}
       \else
           \plainitem{F{\sc IGURE} #1. }#2\par
       \fi
  \endinsert}}
%
%

%
%
\def\Abh#1 {{\sl Abh.\g Math.\g Sem.\g Univ.\g Hamburg\penalty100\ \bf#1\ }}
\def\AMASH#1 {{\sl Acta Math.\g Acad.\g Sci.\g Hung.\penalty100\ \bf#1\ }}
\def\Advances#1 {{\sl Adv.\g Math.\penalty100\ \bf#1\ }}
\def\Annals#1 {{\sl Ann.\g Math.\penalty100\ \bf#1\ }}
\def\AnnComb#1 {{\sl Ann.\g Comb.\penalty100\ \bf#1\ }}
\def\AMM#1 {{\sl Amer.\g Math.\g Monthly\penalty100\ \bf#1\ }}
\def\Archiv#1 {{\sl Arch.}\g {\sl Math.\penalty100\ \bf#1\ }}
\def\ArsComb#1 {{\sl Ars Comb.\penalty100\ \bf#1\ }}
\def\CJM#1 {{\sl Can.\g J.\th Math.\penalty100\ \bf#1\ }}
\def\Comb#1 {{\sl Com\-bi\-na\-to\-ri\-ca\penalty100\ \bf#1\ }}
\def\CPC#1 {{\sl Comb.\g Probab.\g Comput.\penalty100\ \bf#1\ }}
\def\Crelle#1 {{\sl J.}\th {\sl Reine Angew.}\g
    {\sl Math.\penalty100\ \bf#1\ }}
\def\DM#1 {{\sl Discrete Math.\penalty100\ \bf#1\ }}
\def\DAM#1 {{\sl Discrete Appl.\g Math.\penalty100\ \bf#1\ }}
\def\EJC#1 {{\sl Eur.}\g{\sl J.}\g{\sl Comb.\penalty100\ \bf#1\ }}
\def\EJ#1 {{\sl Electronic.}\g{\sl J.}\g{\sl Comb.\penalty100\ \bf#1\ }}
\def\GC#1 {{\sl Graphs Comb.\penalty100\ \bf#1\ }}
\def\IJ#1 {{\sl Isr.\g J.\th Math.\penalty100\ \bf#1\ }}
\def\Inv#1 {{\sl In\-vent.\g math.\penalty100\ \bf#1\ }}
\def\JAlg#1 {{\sl J.}\th {\sl Algorithms\penalty100\ \bf#1\ }}
\def\JCTA#1 {{\sl J.}\th {\sl Comb.}\g {\sl Theory~A\penalty100\ \bf#1\ }}
\def\JCTB#1 {{\sl J.}\th {\sl Comb.}\g {\sl Theory~B\penalty100\ \bf#1\ }}
\def\JGT#1 {{\sl J.}\th {\sl Graph Theory\penalty100\ \bf#1\ }}
\def\BLMS#1 {{\sl Bull.\g Lond.\g Math.\g Soc.\penalty100\ \bf#1\ }}
\def\JLMS#1 {{\sl J.\g Lond.\g Math.\g Soc.\penalty100\ \bf#1\ }}
\def\PLMS#1 {{\sl Proc.\g Lond.\g Math.\g Soc.\penalty100\ \bf#1\ }}
\def\Order#1 {{\sl Order\ \bf#1\ }}
\def\Random#1 {{\sl Random Struct.\g Alg.\penalty100\ \bf#1\ }}
\def\MA#1 {{\sl Math.}\g {\sl Ann.\penalty100\ \bf#1\ }}
\def\MN#1 {{\sl Math.}\g {\sl Nachr.\penalty100\ \bf#1\ }}
\def\MPCPS#1 {{\sl Math.\g Proc.\g Camb.\g Phil.\g Soc.\penalty100\ \bf#1\ }}
\def\MS#1 {{\sl Math.}\g {\sl Scand.\penalty100\ \bf#1\ }}
\def\MZ#1 {{\sl Math.}\g {\sl Zeit.\penalty100\ \bf#1\ }}
\def\BAMS#1 {{\sl Bull.\th Amer.\g Math.\g Soc.\penalty100\ \bf#1\ }}
\def\JAMS#1 {{\sl J.\th Amer.\g Math.\g Soc.\penalty100\ \bf#1\ }}
\def\MAMS#1 {{\sl Mem.\g Amer.\g Math.\g Soc.\penalty100\ \bf#1\ }}
\def\PAMS#1 {{\sl Proc.\g Amer.\g Math.\g Soc.\penalty100\ \bf#1\ }}
\def\SIAM#1 {{\sl SIAM J.\g Discrete Math.\penalty100\ \bf#1\ }}
\def\SLNM#1 {{\sl Springer Lecture Notes in Mathematics\penalty100\ \bf#1\ }}
\def\TAMS#1 {{\sl Trans.\g Amer.\g Math.\g Soc.\penalty100\ \bf#1\ }}
\def\TCSA#1 {{\sl Theor.\g Comput.\g Sci.~A\penalty100\ \bf#1\ }}
%
%
%
%
%
%
%
%
%
%
%
%
%
%
%
\bigskipamount=1\baselineskip plus.3\baselineskip minus.3\baselineskip
\medskipamount=\bigskipamount\divide\medskipamount by 2
\smallskipamount=\medskipamount\divide\smallskipamount by 2 
\medmuskip = 3mu plus 2mu minus 1mu
\thickmuskip = 6mu plus 4mu minus 2mu 
\def\smallbreak{\par \ifdim\lastskip<\smallskipamount
   \removelastskip \penalty-100 \smallskip \fi}
\def\medbreak{\par \ifdim\lastskip<\medskipamount
   \removelastskip \penalty-250 \medskip \fi}
\def\bigbreak{\par \ifdim\lastskip<\bigskipamount
   \removelastskip \penalty-500 \bigskip \fi}
\catcode`@=11        
  \def\raggedbottom{\topskip10pt plus20pt \r@ggedbottomtrue} 
\catcode`@=12        
\def\ge{\geqslant}
\def\le{\leqslant}
\let\elt=\in
\def\in{\mathrel{\mathchoice
   {\raise .7pt \hbox{$\scriptstyle\elt$}}
   {\raise .7pt \hbox{$\scriptstyle\elt$}}
   {\raise .5pt \hbox{$\hskip .5pt\scriptscriptstyle\elt\hskip .5pt$}}
   {\raise.35pt \hbox{$\scriptscriptstyle\elt$}} }}
\let\hasaselt=\owns
\def\owns{\mathrel{\mathchoice
   {\raise .7pt \hbox{$\scriptstyle\hasaselt$}}
   {\raise .7pt \hbox{$\scriptstyle\hasaselt$}}
   {\raise .5pt \hbox{$\hskip .5pt\scriptscriptstyle\hasaselt\hskip .5pt$}}
   {\raise.35pt \hbox{$\scriptscriptstyle\hasaselt$}} }}
\let\exis=\exists
   \def\exists{\exis\>}
\let\nexis=\nexists
   \def\nexists{\nexis\>}
\let\foral=\forall
   \def\forall{\foral\>}
\let\Rightarro=\Rightarrow
   \def\Rightarrow{\>\Rightarro\>}
\let\mi=\min
   \def\min{\mi\>}
\let\ma=\max
   \def\max{\ma\>}
\let\su=\sup
   \def\sup{\su\>}
\let\inff=\inf
   \def\inf{\inff\>}
\mathchardef\to="2221   
\def\proclaim #1.#2 #3\par{\bigbreak
   \noindent{\bf#1.}#2\enspace{\sl#3}\par\bigbreak}
%
\newskip\sectionheadlineskipamount
\sectionheadlineskipamount=8pt plus 2pt minus 1pt
\def\beginsection #1\par{\Bigbreak\centerline{\bold #1}
        \penalty10000\vskip\sectionheadlineskipamount\noindent}
\let\ffootnote=\footnote
\def\footnote#1#2{\ffootnote{#1}{\eightpoint#2\vskip-12pt}}
%
%
\newcount\footnoteno
\def\Footnote#1{{\footnote{${}^{\the\footnoteno}$}{#1}%
   \global\advance\footnoteno by 1}}
\def\item#1#2\par{\parindent=10mm\smallbreak\hang\indent
                  \llap{{\rm #1}\enspace}\ignorespaces#2\par\smallbreak
                  \parindent=7mm}
\def\itemitem#1#2\par{\parindent=10mm\smallbreak
                  \indent\hangindent2\parindent\indent
                  \llap{{\rm #1}\enspace}\ignorespaces#2\par\smallbreak
                  \parindent=7mm}
%
%
%
\pretolerance=0 
\tolerance=2000
\baselineskip=13pt                 
\vsize=200mm                   
\hsize=120mm                   
\hoffset=9mm                   
\parindent=7mm
\relpenalty=2000 \binoppenalty=5000  
\hyphenpenalty=100
\abovedisplayskip=12pt plus3pt minus4pt
\belowdisplayskip=12pt plus3pt minus4pt    
%
%
\belowdisplayshortskip=9pt plus3pt minus3pt
%
%
%
%
%
 \hyphenation{Baum-ord-nung Baum-ord-nun-gen End-ecke End-ecken kur-zen
Kur-zen Graphen-ei-gen-schaft Graphen-ei-gen-schaften he-raus he-raus-ar-bei-ten
he-raus-zu-ar-bei-ten Schnitt-raum}%
 \hyphenation{ac-cess-ible ana-log-ous ana-log-ous-ly ana-lyze ana-lyse
ana-ly-sis answer answers aver-age bundle bundles Buch-ge-sell-schaft col-our
col-ours col-oured col-our-ing col-our-ings con-struct-ible con-struct-ive
con-struct-ive-ly co-rol-lary Co-rol-lary des-cend des-cend-ing Deut-sche
end-li-cher de-fi-ni-tion de-fi-ni-tions De-fi-ni-tion equi-val-ent
equi-val-ence Euler-ian exist-ence every Gra-phen Hamil-ton-ian homeo-mor-phic
homeo-mor-phism homeo-mor-phisms hy-po-thesis hy-po-theses in-ac-cess-ible
ir-regu-lar ir-regu-lar-ity method methods modi-fi-ca-tion mono-chro-matic par-ticu-lar
pro-po-si-tion pro-po-si-tions Pro-po-si-tion regu-lar regu-lar-ity regu-lar-ly
sig-ni-fi-cant sig-ni-fi-cant-ly sig-ni-fi-cance to-po-lo-gical to-po-lo-gical-ly
un-at-tached un-end-li-cher using Using Wis-sen-schaft-li-che}
\userefnumbers\let\newcounter=\newcount
\refno =0
\global \advance \refno by 1 \newcounter \refPrimenden  \refPrimenden =\the \refno 
\global \advance \refno by 1 \newcounter \refCarmesinEF  \refCarmesinEF =\the \refno 
\global \advance \refno by 1 \newcounter \refCWoess  \refCWoess =\the \refno 
\global \advance \refno by 1 \newcounter \refTopSurvey  \refTopSurvey =\the \refno 
\global \advance \refno by 1 \newcounter \refBook  \refBook =\the \refno 
\global \advance \refno by 1 \newcounter \refASS  \refASS =\the \refno 
\global \advance \refno by 1 \newcounter \refTopEnds  \refTopEnds =\the \refno 
\global \advance \refno by 1 \newcounter \refTangleTreeAbstract  \refTangleTreeAbstract =\the \refno 
\global \advance \refno by 1 \newcounter \refDualityGraphsMatroids  \refDualityGraphsMatroids =\the \refno 
\global \advance \refno by 1 \newcounter \refFreudenthal  \refFreudenthal =\the \refno 
\global \advance \refno by 1 \newcounter \refFreudenthalNeu  \refFreudenthalNeu =\the \refno 
\global \advance \refno by 1 \newcounter \refHalinEnds  \refHalinEnds =\the \refno 
\global \advance \refno by 1 \newcounter \refHalinGrid  \refHalinGrid =\the \refno 
\global \advance \refno by 1 \newcounter \refHeuerFullGrid  \refHeuerFullGrid =\the \refno 
\global \advance \refno by 1 \newcounter \refPolat  \refPolat =\the \refno 
\global \advance \refno by 1 \newcounter \refInvLim  \refInvLim =\the \refno 
\global \advance \refno by 1 \newcounter \refGMX  \refGMX =\the \refno 
\global \advance \refno by 1 \newcounter \refRST  \refRST =\the \refno 
\global \advance \refno by 1 \newcounter \refRSTminors  \refRSTminors =\the \refno 
\global \advance \refno by 1 \newcounter \refCTW  \refCTW =\the \refno 
\disablems
\noremarks
\nocomments
\newread\epsffilein    
\newif\ifepsfatend     
\newif\ifepsfbbfound   
\newif\ifepsfdraft     
\newif\ifepsffileok    
\newif\ifepsfframe     
\newif\ifepsfshow      
\epsfshowtrue          
\newif\ifepsfshowfilename 
\newif\ifepsfverbose   
\newdimen\epsfframemargin 
\newdimen\epsfframethickness 
\newdimen\epsfrsize    
\newdimen\epsftmp      
\newdimen\epsftsize    
\newdimen\epsfxsize    
\newdimen\epsfysize    
\newdimen\pspoints     
\pspoints = 1bp        
\epsfxsize = 0pt       
\epsfysize = 0pt       
\epsfframemargin = 0pt 
\epsfframethickness = 0.4pt 
\def\epsfbox#1{\global\def\epsfllx{72}\global\def\epsflly{72}%
   \global\def\epsfurx{540}\global\def\epsfury{720}%
   \def\lbracket{[}\def\testit{#1}\ifx\testit\lbracket
   \let\next=\epsfgetlitbb\else\let\next=\epsfnormal\fi\next{#1}}%
%
%
\def\epsfgetlitbb#1#2 #3 #4 #5]#6{%
   \epsfgrab #2 #3 #4 #5 .\\%
   \epsfsetsize
   \epsfstatus{#6}%
   \epsfsetgraph{#6}%
}%
\def\epsfnormal#1{%
    \epsfgetbb{#1}%
    \epsfsetgraph{#1}%
}%
\newhelp\epsfnoopenhelp{The PostScript image file must be findable by
TeX, i.e., somewhere in the TEXINPUTS (or equivalent) path.}%
\def\epsfgetbb#1{%
%
%
    \openin\epsffilein=#1
    \ifeof\epsffilein
        \errhelp = \epsfnoopenhelp
        \errmessage{Could not open file #1, ignoring it}%
    \else                       
        {
            \chardef\other=12
            \def\do##1{\catcode`##1=\other}%
            \dospecials
            \catcode`\ =10
            \epsffileoktrue         
            \epsfatendfalse     
            \loop               
                \read\epsffilein to \epsffileline
                \ifeof\epsffilein 
                \epsffileokfalse 
            \else                
                \expandafter\epsfaux\epsffileline:. \\%
            \fi
            \ifepsffileok
            \repeat
            \ifepsfbbfound
            \else
                \ifepsfverbose
                    \immediate\write16{No BoundingBox comment found in %
                                    file #1; using defaults}%
                \fi
            \fi
        }
        \closein\epsffilein
    \fi                         
    \epsfsetsize                
    \epsfstatus{#1}%
}%
%
%
\def\epsfclipoff{\def\epsfclipstring{\ifepsfdraft\space clip\fi}}%
\epsfclipoff 
%
%
\def\epsfspecial#1{%
     \epsftmp=10\epsfxsize
     \divide\epsftmp\pspoints
     \ifnum\epsfrsize=0\relax
       \includegraphics{\ifepsfdraft}%
     \else
       \epsfrsize=10\epsfysize
       \divide\epsfrsize\pspoints
       \includegraphics{\ifepsfdraft}%
     \fi
}%
%
\def\epsfframe#1%
{%
  \leavevmode                   
  \setbox0 = \hbox{#1}%
  \dimen0 = \wd0                                
  \advance \dimen0 by 2\epsfframemargin         
  \advance \dimen0 by 2\epsfframethickness      
  \vbox
  {%
    \hrule height \epsfframethickness depth 0pt
    \hbox to \dimen0
    {%
      \hss
      \vrule width \epsfframethickness
      \kern \epsfframemargin
      \vbox {\kern \epsfframemargin \box0 \kern \epsfframemargin }%
      \kern \epsfframemargin
      \vrule width \epsfframethickness
      \hss
    }
    \hrule height 0pt depth \epsfframethickness
  }
}%
\def\epsfsetgraph#1%
{%
   %
   %
   \leavevmode
   \hbox{
     \ifepsfframe\expandafter\epsfframe\fi
     {\vbox to\epsfysize
     {%
        \ifepsfshow
            \vfil
            \hbox to \epsfxsize{\epsfspecial{#1}\hfil}%
        \else
            \vfil
            \hbox to\epsfxsize{%
               \hss
               \ifepsfshowfilename
               {%
                  \epsfframemargin=3pt 
                  \epsfframe{{\tt #1}}%
               }%
               \fi
               \hss
            }%
            \vfil
        \fi
     }%
   }}%
   %
   %
   \global\epsfxsize=0pt
   \global\epsfysize=0pt
}%
%
%
\def\epsfsetsize
{%
   \epsfrsize=\epsfury\pspoints
   \advance\epsfrsize by-\epsflly\pspoints
   \epsftsize=\epsfurx\pspoints
   \advance\epsftsize by-\epsfllx\pspoints
%
%
   \epsfxsize=\epsfsize{\epsftsize}{\epsfrsize}%
   \ifnum \epsfxsize=0
      \ifnum \epsfysize=0
	\epsfxsize=\epsftsize
        \epsfysize=\epsfrsize
	\epsfrsize=0pt
%
%
      \else
	\epsftmp=\epsftsize \divide\epsftmp\epsfrsize
	\epsfxsize=\epsfysize \multiply\epsfxsize\epsftmp
	\multiply\epsftmp\epsfrsize \advance\epsftsize-\epsftmp
	\epsftmp=\epsfysize
	\loop \advance\epsftsize\epsftsize \divide\epsftmp 2
	\ifnum \epsftmp>0
	   \ifnum \epsftsize<\epsfrsize
           \else
	      \advance\epsftsize-\epsfrsize \advance\epsfxsize\epsftmp
           \fi
	\repeat
	\epsfrsize=0pt
      \fi
   \else
     \ifnum \epsfysize=0
       \epsftmp=\epsfrsize \divide\epsftmp\epsftsize
       \epsfysize=\epsfxsize \multiply\epsfysize\epsftmp
       \multiply\epsftmp\epsftsize \advance\epsfrsize-\epsftmp
       \epsftmp=\epsfxsize
       \loop \advance\epsfrsize\epsfrsize \divide\epsftmp 2
       \ifnum \epsftmp>0
	  \ifnum \epsfrsize<\epsftsize
          \else
	     \advance\epsfrsize-\epsftsize \advance\epsfysize\epsftmp
          \fi
       \repeat
       \epsfrsize=0pt
     \else
       \epsfrsize=\epsfysize
     \fi
   \fi
}%
%
%
\def\epsfstatus#1{
   \ifepsfverbose
     \immediate\write16{#1: BoundingBox:
                  llx = \epsfllx\space lly = \epsflly\space
                  urx = \epsfurx\space ury = \epsfury\space}%
     \immediate\write16{#1: scaled width = \the\epsfxsize\space
                  scaled height = \the\epsfysize}%
   \fi
}%
%
%
{\catcode`\%=12 \global\let\epsfpercent=
\global\def\epsfatend{(atend)}%
%
%
%
%
%
%
%
\long\def\epsfaux#1#2:#3\\%
{%
   \def\testit{#2}
   \ifx#1\epsfpercent           
       \ifx\testit\epsfbblit    
            \epsfgrab #3 . . . \\%
            \ifx\epsfllx\epsfatend 
                \global\epsfatendtrue
            \else               
                \ifepsfatend    
                \else           
                    \epsffileokfalse
                \fi
                \global\epsfbbfoundtrue
            \fi
       \fi
   \fi
}%
%
%
\def\epsfempty{}%
\def\epsfgrab #1 #2 #3 #4 #5\\{%
   \global\def\epsfllx{#1}\ifx\epsfllx\epsfempty
      \epsfgrab #2 #3 #4 #5 .\\\else
   \global\def\epsflly{#2}%
   \global\def\epsfurx{#3}\global\def\epsfury{#4}\fi
}%
%
%
\def\epsfsize#1#2{\epsfxsize}%
%
%


\lineskiplimit=-3pt
\def\varprojlim{\displaystyle{\lim_{\vbox{\hbox{$\longleftarrow$}}}\,}}
\hbox{}\vskip2pt\centerline{\bigbold Ends and tangles}
\smallskip\centerline{\it To the memory of Rudolf Halin}
\vskip 4mm
\centerline{Reinhard Diestel}
\disablems
\noremarks
\nocomments

\def\A{{\aleph_0}}
\def\AT{$\A$-tangle}
\def\C{{\cal C}}
\def\D{{\cal D}}
\def\F{{\cal F}}
\def\O{{\cal O}}
\def\S{{\cal S}}
\def\T{{\cal T}}
\def\U{{\cal U}}
\def\X{{\cal X}}

\def\lowfwd #1#2#3{{\mathop{\kern0pt #1}\limits^{\kern#2pt\raise.#3ex \vbox to 0pt{\hbox{$\scriptscriptstyle\rightarrow$}\vss}}}}

\def\lowbkwd #1#2#3{{\mathop{\kern0pt #1}\limits^{\kern#2pt\raise.#3ex
\vbox to 0pt{\hbox{$\scriptscriptstyle\leftarrow$}\vss}}}}
\def\vS{\vec S}
\def\vs{\lowfwd s{1.5}1}
\def\sv{\lowbkwd s{1.5}1}

\def\uf{ultrafilter}
\def\restricts{\!\restriction\!}
\def\fXX{f_{X'\!,X}}
\def\UX{(\,\U_X\mid X\in\X\,)}

\def\xxxTop{1}
\def\xxxUFs{2}
\def\xxxTau{3}
\def\xxxBlocks{4}

\def\secTangles{1}
\def\secUFs{2}
\def\secUFTs{3}
\def\secTop{4}
\def\secBlocks{5}
\def\secOutlook{6}

\def\xxxSmall{\secTangles.1}
\def\xxxP{\secTangles.2}
\def\xxxFactfive{\secTangles.3}
\def\xxxFactthree{\secTangles.4}
\def\xxxFactfour{\secTangles.5}
\def\xxxDirections{\secTangles.6}
\def\xxxEndTangles{\secTangles.7}
\def\xxxCombEx{\secTangles.8}
\def\xxxStarEx{\secTangles.9}
\def\xxxSymDiff{\secTangles.10}
\def\xxxInfiniteB{\secTangles.11}

\def\xxxInducedUFs{\secUFs.1}
\def\xxxCommute{\secUFs.2}
\def\xxxAllTangles{\secUFs.3}
\def\xxxPrincipal{\secUFs.4}

\def\xxxg{\secUFTs.1}
\def\xxxgf{\secUFTs.2}
\def\xxxX{\secUFTs.3}
\def\xxxDirect{\secUFTs.4}
\def\xxxtau{\secUFTs.5}
\def\xxxXtau{\secUFTs.6}
\def\xxxSurjective{\secUFTs.7}

\def\xxxcts{\secTop.1}
\def\xxxTanglesExist{\secTop.2}
\def\xxxTheta{\secTop.3}
\def\xxxBasis{\secTop.4}

\def\xxxRay{\secBlocks.1}
\def\xxxUFnotClosed{\secBlocks.2}
\def\xxxComplete{\secBlocks.3}
\def\xxxTK{\secBlocks.4}

\def\eqomegatau{{\rm(1)}}
\def\eqomegaf{{\rm(2)}}
\def\eqtauf{{\rm(3)}}
\def\eqtauUtauX{{\rm(4)}}
\def\eqtauy{{\rm(5)}}
\def\eqomegay{{\rm(6)}}
\def\eqbasic{{\rm(7)}}
\def\eqopensets{{\rm(8)}}
\def\equps{{\rm(9)}}
\def\eqK{{\rm(10)}}

\def\figDiagram{1}
\def\figComb{2}
\def\figStar{3}
\def\figUFs{4}
\def\figBigDiagram{5}

\bigskip\medskip
{\narrower\narrower\ninepoint\noindent
   We show that an arbitrary infinite graph can be compactified by its \AT s in much the same way as the ends of a locally finite graph compactify it in its Freudenthal compactification. In general, the ends then appear as a subset of its \AT s.

The \AT s of a graph are shown to form an inverse limit of the ultrafilters on the sets of components obtained by deleting a finite set of vertices. The \AT s that are ends are precisely the limits of principal ultrafilters.

The \AT s that correspond to a highly connected part, or $\A$-block, of the graph are shown to be precisely those that are closed in the topological space of its finite-order separations.\par}

\beginsection Introduction

Much of Halin's legacy in graph theory stems from the fact that, in his seminal paper~[\the\refHalinEnds] of~1964, he initiated the study of {\it ends\/} for infinite graphs. Our aim in this paper is to unify this notion with that of a {\it tangle\/} introduced by Robertson and Seymour~[\the\refGMX] in~1991. It turns out that Halin's ends can be viewed as a special case of tangles of infinite order. These can, in turn, be used to compactify an arbitrary infinite graph in the same way as ends compactify locally finite graphs in their well-known Freudenthal compactification.

Inspired by Carath\'eodory's notion of {\it Primenden\/} of regions in the complex plane~[\the\refPrimenden], but unaware of Freudenthal's~[\the\refFreudenthal] generalization of these to more general locally compact Hausdorff spaces, Halin~[\the\refHalinEnds] defined an {\it end\/} of a graph~$G$ as an equivalence class of {\it rays\/}, or 1-way infinite paths, in~$G$. Here, two rays are {\it equivalent\/} if no finite set of vertices separates them in~$G$.

Halin does not require these graphs to be locally finite. If they are, his notion of an end is equivalent to Freudenthal's. If a graph is locally finite and connected, there is a natural topology that makes it and its ends into a compact space, its {\it Freudenthal compactification\/}~[\the\refTopSurvey,\th \the\refBook]. The rays of an end, in Halin's definition, then converge to their ends. If the graph is not locally finite, Halin's notion of an end no longer agrees with Freudenthal's~[\the\refFreudenthalNeu]%
   \COMMENT{}
   but is more general~[\the\refTopEnds]. There is no longer an obvious topology on the graph and its ends, in either definition, that makes rays converge to `their' ends~-- let alone one that makes the graph and its ends compact.

An end $\omega$ of a graph~$G$ orients its separations $\{A,B\}$ of finite order in that either every ray from~$\omega$ has a tail in~$A$, or every ray from~$\omega$ has a tail in~$B$. These orientations of finite-order separations are {\it consistent\/} in a number of ways; for example, if $\{C,D\}$ is another such separation with $C\sub A$ and $D\supe B$, then $\omega$ must orient $\{C,D\}$ towards~$D$ if it orients $\{A,B\}$ towards~$B$.

Robertson and Seymour~[\the\refGMX], independently, introduced another notion of consistently orienting all the low-order separations of~$G$: the notion of a {\it tangle\/}. It is easy to see that every end~$\omega$ defines an \AT, one that orients all the separations of finite order: if we orient them all towards the side where the rays in~$\omega$ have their tails, we obtain an \AT. Conversely, every \AT\ of a locally finite and connected graph~$G$ is defined by an end in this way. Thus, if $G$ is locally finite, its \AT s are just another way of identifying the points at infinity in its Freudenthal compactification~$|G|$.

When $G$ is not locally finite, however, things get interesting. Now adding its ends no longer compactifies~$G$.%
   \COMMENT{}
   But there can also be \AT s that are not defined by an end. And it turns out that adding these, in addition to the ends, does compactify~$G$. More precisely, we shall define a topology on the union of~$G$, viewed as a 1-complex, and its set $\Theta$ of \AT s so that the resulting space $|G| = G\cup\Theta$ satisfies the following:

\proclaim Theorem~\xxxTop. Let $G$ be any graph.
\pitem{i} $|G|$ is a compact space in which $G$ is dense and $|G|\sm G$ is totally disconnected.
\pitem{ii} If $G$ is locally finite and connected, then all its \AT s are ends, and $|G|$ coincides with the Freudenthal compactification of~$G$.
   \enditem

\noindent
   This compactification of $G$ differs from its \SC\ compactification%
   \COMMENT{}
   and from other compactifications that have been suggested for graphs that are not locally finite, e.g. by Cartwright, Soardi and Woess~[\the\refCWoess] or by Polat~[\the\refPolat].

In order to prove Theorem~\xxxTop, we shall need to understand the tangles that are not ends. When $X$ is a finite set of vertices of~$G$, then every bipartition of the set $\C_X$ of the components of $G-X$ defines a finite-order separation of~$G$: if $\C_X = \C_1\cup \C_2$, say, then $\{\bigcup\C_1\cup X, \bigcup\C_2\cup X\}$ is such a separation. Every ultrafilter $U$ on $\C_X$ contains exactly one $\C_i$ from such a bipartition. If we think of $U$ as orienting the corresponding separation of $G$ towards this~$\C_i$, and there are no other finite-order separations in~$G$, then $U$ defines an \AT\ of~$G$. Conversely, every \AT\ of~$G$ will orient all the separations of the above form in such a way that the $\C_i$ it points to form an \uf\ on~$\C_X$. We shall call all the \uf s on such sets~$\C_X$ with $X$ finite the {\it ultrafilters of cofinite components\/} of~$G$.

Of course, there will normally be other finite-order separations $\{A,B\}$ of~$G$.\penalty-200\ But each of these also has the above form, with $X:= A\cap B$. And it turns out that the \uf s $U_X$ which a tangle defines for different choices of~$X$ are all compatible in a simple and natural way: they form the limits $(\,U_X\mid X\in\X\,)$ of a natural inverse system $(\,\U_X\mid X\in\X\,)$ of the sets of ultrafilters of cofinite\penalty-200\ components (where $\X$ is the set of finite sets of vertices of~$G$). And conversely, every such limit of \uf s comes from a tangle:

\proclaim Theorem~\xxxUFs. Let $G$ be any graph.
\pitem{i} The \AT s of $G$ are precisely the limits of the inverse system of its sets of \uf s of cofinite components. 
\pitem{ii} The ends of $G$ are precisely those of its \AT s whose \uf s of cofinite components are all principal.
 \enditem

When a tangle is defined by an end, every ray in that end can be used as an `oracle' to determine which of the two orientations of a given separation lies in the tangle: it is the orientation that points to a tail of that ray. For tangles that are not defined by an end we can use an \uf\ in a similar way:

\proclaim Theorem \xxxTau.
Every \AT\ $\tau$ in~$G$ satisfies exactly one of the following:
 \medskip
 \plainitem{$\bullet$} There is a ray $R$ in $G$ such that every separation in~$\tau$ points to a tail of~$R$.
 \smallskip
 \plainitem{$\bullet$} There is a non-principal \uf~$U$ of cofinite components in~$G$ such that every separation in~$\tau$ points to an element of~$U$.
 \enditem

In a finite graph, the tangles of order some fixed~$k\in\N$, those that orient every separation of order~$<k$, can be thought of as pointing towards some `highly connected substructure' of that graph. This is clearly not the case for all \AT s: the ends of a tree, for example, are hardly highly connected substructures. However, we can identify those \AT s of an infinite graph~$G$ that do point to a highly connected substructure in a meaningful sense: they are the \AT s that are closed in the set~$\vS$ of all oriented separations of finite order of the graph, with respect to a natural topology on this set.

End tangles of trees are not closed in this topology, but \AT s pointing to some fixed infinite complete subgraph, for example, are. More generally, tangles pointing to a fixed $\A$-block are closed: a {\it $\kappa$-block\/} in a graph~$G$ is a maximal set of at least $\kappa$ vertices no two of which can be separated in~$G$ by fewer than $\kappa$ vertices. As it turns out, the ends defined by an $\A$-block are precisely the \AT s that are closed in~$\vS$:

\proclaim Theorem~\xxxBlocks.
Let $G$ be any graph.
\pitem{i} $\!$The \AT s of $G$ that are not ends are never closed in~$\vS$.
\pitem{ii} $\!$An end tangle of $G$ is closed in~$\vS$ if and only if it is defined by an $\A$-block.

The paper is organized as follows. We begin in Section~\secTangles\ with a review of tangles, especially those of infinite order, and their relationship to ends. In Section~\secUFs\ we introduce the inverse system of the sets $\U_X$ of ultrafilters of cofinite components, and prove Theorem~\xxxUFs. In Section~\secUFTs\ we take a closer look at how tangles not defined by an end arise in this inverse system: we show that each of them is already determined by a single non-principal \uf\ among those of which it is a limit. In Section~\secTop\ we topologize the set of \AT s by viewing it as~$\varprojlim\U_X$, with the $\U_X$ carrying the \SC\ topology. We then extend this space to include $G$ itself, and prove Theorem~\xxxTop. In Section~\secBlocks, finally, we introduce a topology on the set~$\vS$ of finite-order separations of a graph, and prove Theorem~\xxxBlocks. We close with a short section on the potential for applications of Theorem~\xxxTop\ and possible further questions, one of which appears to be quite far-reaching.

Any graph-theoretic notation not explained here can be found in~[\the\refBook]. Inverse systems and inverse limits, including their topology, are explained in~[\the\refInvLim].

\beginsection \secTangles. Tangles and ends

A {\it separation\/} of a graph $G=(V,E)$ is a set $\{A,B\}$ such that $A\cup B=V$ and $G$ has no edge between $A\sm B$ and $B\sm A$. The {\it order\/} of a separation $\{A,B\}$, and of its orientations (see below), is the cardinal number~$|A\cap B|$.

The ordered pairs $(A,B)$ and $(B,A)$ are the {\it orientations\/} of a separation $\{A,B\}$ and are also called {\it (oriented) separations\/}. Given a set $S$ of separations, we write $\vS$ for the set of their orientations. An {\it orientation of~$S$} is a subset $O$ of $\vS$ that contains for every $\{A,B\}\in S$ exactly one of $(A,B)$ and~$(B,A)$. 

Mapping every $(A,B)\in\vS$ to its {\it inverse\/} $(B,A)$ is an involution on~$\vS$ that reverses the partial ordering
 $$(A,B)\le (C,D) :\Leftrightarrow A\subseteq C \hbox{ and } B\supseteq D,$$
 since the above is equivalent to $(D,C)\le (B,A)$. Informally, we think of $(A,B)$ as {\it pointing towards}~$B$ and {\it away from}~$A$.
   Similarly, if $(A,B)\le (C,D)$, then $(A,B)$ {\it points towards} $\{C,D\}$ and its orientations, while $(C,D)$ {\it points away from} $\{A,B\}$ and its orientations.

A~set $\sigma$ of oriented separations is a {\it star\/} if they all point towards each other: if $(A,B)\le (B',A')$ for all distinct $(A,B), (A',B')\in\sigma$. Note that if $\sigma$ is a star then $\bigcap\,\{\,B\mid (A,B)\in\sigma\,\}$ contains $A'\cap B'$ for every $(A',B')\in\sigma$.%
   \COMMENT{}

A~set of oriented separations is {\it consistent\/} if no two of them point away from each other: if it contains no distinct separations $(B,A)$ and $(C,D)$ with $(A,B) < (C,D)$. For an orientation $O$ of~$S$, consistency is tantamount to being closed down in~$\vS\,$: that $(A,B) < (C,D)\in O$ with $\{A,B\}\in S$ implies $(A,B)\in O$.

For the rest of this paper, let $G=(V,E)$ be a fixed infinite graph. Let $S=S_\A$ be the set of all its separations of finite order,%
   \COMMENT{}
   and let $\S$ denote the set of stars in~$\vS$. Let $\X$ be the set of finite subsets of~$V$. As usual, we write $\Omega = \Omega(G)$ for the set of ends of~$G$, defined as in the Introduction or in~[\the\refBook].

\proclaim Lemma \xxxSmall.
Every consistent orientation $O$ of $S$ contains every separation $(A,V)$ of $G$ with $A$ finite.

\proof
Pick $v\in V\sm A$, and let $A' := A\cup\{v\}$. Since $O$ contains one of~$(A',V)$ and~$(V,A')$, and $(A,V) < (A',V)$ as well as $(A,V) < (V,A')$, the consistency of $O$ requires that $(V,A)\notin O$. Hence $(A,V)\in O$, as claimed.
 \endproof

  Given a set $\F\sub\S$ of stars of separations, we call a consistent orientation of~$S$ an {\it $\F$-tangle\/} if it has no subset in~$\F$.%
   \COMMENT{}
Let us consider some particular choices for~$\F$. For integers~$n\ge 1$, let
 $$\T_n := \big\{\,\{(A_1,B_1),\ldots,(A_n,B_n)\}\in\S: B_1\cap\ldots\cap B_n \hbox{ is finite}\,\big\}\rlap.$$
   These $(A_i,B_i)$ need not be distinct, so $\T_n$ is a set of stars in~$\vS$ of up to~$n$ separations each. In particular, $\T_m\sub\T_n$ for all $m\le n$, so every $\T_n$-tangle of~$S$%
   \COMMENT{}
    is also a $\T_m$-tangle.

Let us note the following observation for later use:

\proclaim Lemma~\xxxP.
Any $\T_3$-tangle $O$ of~$S$ containing two separations $(A,B), (A',B')$ also contains $(A\cup A', B\cap B')$.

\proof
$\{A\cup A', B\cap B'\}$ is a clearly a separation of~$G$. It lies in~$S$, because its separator is a subset of $(A\cap B)\cup (A'\cap B')$ and hence has finite order.  Suppose $(B\cap B', A\cup A')\in O$. Since also $(A\cap B', B\cup A')\le (A,B)\in O$ lies in~$O$, because $O$ is consistent, and $(A',B')$ does by assumption, we have found three separations in~$O$ forming a star in~$\T_3$, a contradiction.
   \endproof

Robertson and Seymour~[\the\refGMX] defined tangles slightly differently, as follows.%
   \Footnote{Formally, their definition of a $k$-tangle (for finite $k$ and~$G$) differs slightly from our definition below, but is easily seen to be equivalent.}
   Given a cardinal~$\kappa$, a {\it tangle of order~$\kappa$}, or {\it $\kappa$-tangle\/}, of $G$ is an orientation of the set $S_\kappa$ of its separations of order $<\kappa$ that has no subset of the form
 $$\{\,(A_1,B_1), (A_2,B_2), (A_3,B_3) : G[A_1]\cup G[A_2]\cup G[A_3] = G\,\},\eqno(\rm T)$$
 where $G[A_i]$ denotes the subgraph of~$G$ induced by~$A_i$. As before, the $(A_i,B_i)$ need not be distinct, so in particular every $\kappa$-tangle is consistent.%
   \COMMENT{}

Our $\T_3$-tangles of~$S=S_\A$, however, coincide with its $\A$-tangles as defined by Robertson and Seymour:

\proclaim Lemma~\xxxFactfive.
The $\T_3$-tangles of $S$ are precisely the $\A$-tangles of~$G$.

\proof
As remarked, both $\T_3$-tangles of $S$ and $\A$-tangles of $G$ are consistent orientations of~$S$. It was shown in~[\the\refDualityGraphsMatroids, Lemma~4.2] that consistent orientations of $S$ with no star subset as in~(T) have no subset as in~(T) at all, star or not. Since any set as in~(T) satisfies 
 $$\textstyle\bigcap_{i=1}^3 B_i = V\cap \bigcap_{i=1}^3 B_i = \bigcup_{i=1}^3 A_i \cap \bigcap_{i=1}^3 B_i\sub \bigcup_{i=1}^3 (A_i\cap B_i),$$%
   \COMMENT{}
 which is finite, all $\T_3$-tangles of $S$ are $\A$-tangles of~$G$.

Conversely, suppose some consistent orientation $O$ of~$S$ is not a $\T_3$-tangle but contains a star $\sigma = \{(A_1,B_1), (A_2,B_2), (A_3,B_3)\}\in\T_3$. We show that $O$ also has a subset as in~(T), and hence is not an $\A$-tangle. As $X\! := B_1\cap B_2\cap B_3$ is finite,%
   \COMMENT{}
   we have $(A',B_1) := (A_1\cup X, B_1)\in\vS$.%
   \COMMENT{}
   Using that $\sigma$ is a star, one easily checks that $G[A']\cup G[A_2]\cup G[A_3] = G$. Hence if $(A',B_1)\in O$, then
 $$\{(A',B_1), (A_2,B_2), (A_3, B_3)\}$$
is our desired subset of $O$ as in~(T). But if not, then $(B_1,A')\in O$. But then $\{(A_1,B_1), (B_1,A')\}$ is a subset of $O$ as in~(T), since $G[A_1]\cup G[B_1] = G$.
 \endproof

We noted earlier that, trivially, every $\T_n$-tangle of~$S$ is also a $\T_m$-tangle for all $m\le n$. By the particular nature of our $S=S_\A$, we also have a converse:

\proclaim Lemma~\xxxFactthree.
Every $\T_n$-tangle of~$S$ is also a $\T_{n+1}$-tangle, for all $n\ge 3$.

\proof
Let $O\sub\vS$ be a $\T_n$-tangle of~$S$, and let $\sigma\in\T_{n+1}$ be given; we have to show that $\sigma\not\sub O$. Pick distinct $(A_1,B_1), (A_2,B_2)\in\sigma$.%
   \COMMENT{}
   We may assume that $(A_1,B_1), (A_2,B_2)\in O$, as otherwise $\sigma\not\sub O$ as desired.

Let $A:= A_1\cup A_2$ and $B:= B_1\cap B_2$, and let $\sigma'$ be obtained from~$\sigma$ by replacing $(A_1,B_1)$ and $(A_2,B_2)$ with~$(A,B)$. It is easy to check that, since $\sigma$ is a star in~$\S$, so is~$\sigma'$.%
   \COMMENT{}
   Hence $\sigma'\in\T_n$, giving $\sigma'\not\sub O$ by the choice of~$O$.

Hence if $\sigma\sub O$%
   \COMMENT{}
   then $(A,B)\notin O$, and thus $(B,A)\in O$. (Here we use that $\{A,B\}\in S$, which can fail for arbitrary~$S$ but clearly holds for our $S=S_\A$.) But now $O$ contains the set $\{(B,A), (A_1,B_2), (A_2, B_2)\}\in\T_3\sub\T_n$, contrary to its definition. Thus, $\sigma\not\sub O$ as desired.
   \endproof

For $\T_{<\aleph_0}:= \bigcup_{n=3}^\infty \T_n$, Lemma~\xxxFactthree\ implies

\proclaim Corollary \xxxFactfour.
The $\T_{<\aleph_0}$-tangles of $S$ are exactly its $\T_3$-tangles.\noproof

As we have seen, all choices of $\T_n$ or of~$\T_{<\aleph_0}$ as $\F$ yield the same $\F$-tangles: if a consistent orientation of $S$ avoids all 3-stars whose target sides have a finite intersection, it avoids all finite such stars. In view of Lemma~\xxxFactfive, we shall from now on refer to all these $\F$-tangles of~$S$ as the {\it $\A$-tangles\/} in~$G$,%
   \COMMENT{}
   and write $\Theta = \Theta(G)$ for the set of all these.

But how about excluding infinite stars of separations as well? Let
 $$\T := \big\{\,\{(A_i,B_i)\mid i\in I\}\in\S: \bigcap_{i\in I} B_i \hbox{ is finite}\,\big\},$$
where the $I$ are arbitrary index sets. Thus $\T_{<\aleph_0}\!\sub\T$, and hence all $\T$-tangles are \AT s. Unlike in Lemma~\xxxFactthree, the converse does not hold. But before we give an example, let us show something more surprising: the $\T$-tangles of $S$ correspond precisely to the ends of~$G$!

One way of this correspondence is easy. Given an end~$\omega$ of~$G$, exactly one of the two orientations $(A,B)$ of each separation in~$S$ has the property that $B$ contains a tail of every ray in~$\omega$ (for which we say  that $\omega$ {\it lives in\/}~$B$): this is immediate from the notion of ray-equivalence in the definition of an end. Hence $\omega$ defines an orientation of~$S$, which is easily seen to be consistent. This orientation~$\tau$ has no subset $\sigma = \{(A_i,B_i)\mid i\in I\}$ in~$\T$. Indeed, as $X:= \bigcap_i B_i$ is finite, our end $\omega$ lives in some component $C$ of $G-X$.%
   \COMMENT{}
   If $C\cap A_i\ne\es$, say,%
   \Footnote{To simplify notation, we do not always distinguish between graphs and their vertex sets.\looseness=-1} 
   then $C\sub A_i\sm B_i$, because $A_i\cap B_i\sub X$ since $\sigma$ is a star.%
   \COMMENT{}
   Therefore $(B_i,A_i)\in\tau$, and hence $\sigma\not\sub\tau$ as claimed. 

We have thus defined a map
 $$\omega\mapsto\tau_\omega\eqno\eqomegatau$$
from the ends of~$G$ to its $\A$-tangles, whose images are in fact $\T$-tangles. We shall say that $\tau_\omega$ is {\it defined by\/}~$\omega$.

The map in~\eqomegatau\ is clearly injective. Indeed, distinct ends are, by definition, distinguished by some $X\in\X$ in the sense that they live in different components of~$G-X$. Then any separation $\{A,B\}\in S$ with $A\cap B = X$ for which these components lie on different sides will get oriented differently by the corresponding two \AT s.

Conversely,%
   \COMMENT{}
   every $\T$-tangle is defined by an end in this way. To prove this we need a result from~[\the\refTopEnds] (see also~[\the\refRST]), which requires another definition.

A {\it direction\/} in $G$ is a function $f$ that assigns to every finite set $X\sub V$ one of the components of $G-X$ so that $f(X')\sub f(X)$ whenever $X\sub X'$. Clearly, every end $\omega$ of $G$ defines a direction~$f$ by taking as $f(X)$ the unique component of $G-X$ in which $\omega$ lives. It was shown in~[\the\refTopEnds] that this map
 $$\omega\mapsto f_\omega\eqno\eqomegaf$$
   is a bijection: not only do different ends define different directions (which is immediate), but every direction is defined by an end in the way indicated. Hence all we have to show is that the $\T$-tangles in a graph correspond to its directions:

\proclaim Lemma~\xxxDirections.
For every $\T$-tangle $\tau$ of~$S$ there is a unique direction $f = f_\tau$ in~$G$ such that, for every $X\in\X$ and every component $C$ of $G-X$, we have $(V\sm C, X\cup C)\in\tau$ if and only if $C = f(X)$.
  This map
 $$\tau\mapsto f_\tau\eqno\eqtauf$$
 is a bijection from the $\T$-tangles of~$S$ to the directions in~$G$, which commutes with the maps from {\rm\eqomegatau} and~{\rm\eqomegaf}.

\proof
To define $f$ given~$\tau$, let $X\in\X$ be given. There is a unique component $C$ of $G-X$ such that $(V\sm C, X\cup C)\in\tau$: existence follows from $\tau\in\T$,%
   \COMMENT{}
   while uniqueness follows from the consistency of~$\tau$. Let $f(X) := C$.

\goodbreak

To show that $f\: X\mapsto C$ is a direction of~$G$, consider $X\sub X'\in\X$. As $C':= f(X')$ is connected, it lies inside a component of~$G-X$. If this component was not~$C$, the separations $(V\sm C, X\cup C), (V\sm C', X'\cup C')\in\tau$ would make $\tau$ inconsistent, contradicting our assumptions.

 \figure \figDiagram.
 Ends, tangles, and directions
 (Diagram; 1000)

Clearly, our map $\tau\mapsto f$ composes with the map $\omega\mapsto\tau$ from~\eqomegatau\ as shown in Figure~\figDiagram, to yield the bijection $\omega\mapsto f$ from~\eqomegaf. In particular, it is surjective.%
   \COMMENT{}
   In order to show that it is injective,%
   \COMMENT{}
   consider distinct $\T$-tangles of~$S$, say $\tau$ and~$\tau'$. Let $\{A,B\}\in S$ be a separation which they orient differently, with $(A,B)\in\tau$ and $(B,A)\in\tau'$ say. Since $\tau$ is consistent, $f_\tau$ maps $X= A\cap B$ to a component of $G-X$ contained in~$B\sm A$.%
   \COMMENT{}
   Similarly, $f_{\tau'}$ maps $X$ to a component of $G-X$ contained in~$A\sm B$. Thus $f_\tau\ne f_{\tau'}$, as desired.
 \endproof

\proclaim Corollary \xxxEndTangles.
The map $\omega\to\tau_\omega$ defined in~\eqomegatau\ is a bijection from the ends of~$G$ to its $\T$-tangles.\noproof

Corollary~\xxxEndTangles\ says that those $\T_{<\aleph_0}$-tangles of~$S$ that are even $\T$-tangles are precisely the \AT s that are defined by an end via~\eqomegatau. We shall call these \AT s the {\it end tangles\/} of~$G$. Let us look at an example.

 \figure \figComb.
 All $\A$-tangles in this graph are end tangles
 (Comb; 1000)

\proclaim Example \xxxCombEx.
{ \rm
If $G$ is the {\it comb of rays\/} shown in Figure~\figComb, then every \AT\ $\tau$ in $G$ is an end tangle. Indeed, if there exists an $i\in\N$ such that $\omega_i$ lives in~$B$ for every $(A,B)\in\tau$, then $\omega_i\mapsto\tau$ in the map of~\eqomegatau, so $\tau$ is an end tangle. Suppose, then, that there is no such~$i$.
   \endgraf
We claim that $\omega$ lives in every $B$ with $(A,B)\in\tau$, so that $\omega\mapsto\tau$ in~\eqomegatau, again making $\tau$ an end tangle. If not then, for some $(A,B)\in\tau$, the end $\omega$ lives in~$A$. Then only finitely many~$\omega_i$ live in~$B$. By assumption,%
   \COMMENT{}
   each of these~$\omega_i$ lives in~$A_i$ for some $(A_i,B_i)\in\tau$. By the consistency of~$\tau$, these $(A_i,B_i)$ can be chosen so as to form a (finite) star $\sigma\sub\tau$ together with~$(A,B)$.%
   \COMMENT{}
   But the intersection of $B$ with these~$B_i$ is finite. Hence $\tau\supe\sigma\in\T_{<\aleph_0}$, so $\tau$ is not an \AT, contrary to our assumption.%
   \COMMENT{}
   \endproof
   }

Next, let us see an example of an \AT\ that is not an end tangle: a consistent orientation of~$S$ that has infinite but no finite stars in~$\T$.

 \figure \figStar.
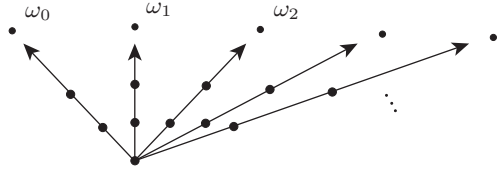
 This graph has both end and ultrafilter tangles
 (Star; 1000)

\proclaim Example \xxxStarEx.
{\rm
If $G$ is the graph shown in Figure~\figStar, it has \AT s that are not end tangles. Indeed, let $U$ be any \uf\ on~$\N$. Every separation $\{A,B\}\in S$ induces a bipartition of~$\N$ into
 $$\bar A = \{\,i\in\N\mid\omega_i\hbox{ lives in }A\,\}\quad\hbox{and}\quad
   \bar B = \{\,i\in\N\mid\omega_i\hbox{ lives in }B\,\}.$$
 As $U$ is an \uf, exactly one of these sets is an element of~$U$. Hence 
 $$\tau = \{\,(A,B)\in\vS\mid \bar B\in U\,\}$$
 is an orientation of~$S$. Since the intersection of two sets in~$U$ lies in $U$ and hence is non-empty, $\tau$~is consistent.%
   \COMMENT{}
   Similarly, let $\sigma\sub\tau$ be a finite star. The set of all $i\in\N$ whose $\omega_i$ lives in every $B$ with $(A,B)\in\sigma$ is a finite intersection of sets in~$U$, and hence is non-empty. Consider any $i$ in this set, and a ray in~$\omega_i$. This ray has a tail outside every~$A$ with $(A,B)\in\sigma$, and hence has a tail in $\bigcap\,\{B\mid (A,B)\in\sigma\}$. In particular, this intersection is infinite, and hence $\sigma\notin\T_{<\aleph_0}$. Thus, $\tau$~is indeed an \AT.
   \endgraf
  If $U$ is a principal \uf\ generated by~$\{n\}$, say, then $\tau$ is an end tangle defined by~$\omega_n$. If $U$ is a non-principal ultrafilter, then $\tau$ is not defined by any~$\omega_i$ and hence is not an end tangle.
   \endproof
   }

We shall see in Section~\secUFTs\ that every \AT\ that is not an end tangle is defined by a non-principal ultrafilter in a way similar to Example~\xxxStarEx.

\medbreak

We conclude this section with a couple of simple lemmas about \AT s. The first is that if we change a separation in an \AT~$\tau$ only finitely, the resulting separation will again lie in~$\tau$. For sets $A,A'\sub V$ let us write $A\sim A'$ if their symmetric difference is finite.

\proclaim Lemma \xxxSymDiff.
Let $\tau$ be an \AT\ of~$S$ and $(A,B)\in\tau$. Let $(A',B')$ be a separation such that $A'\sim A$ and $B'\sim B$. Then $(A',B')\in\tau$.

\proof
It suffices to show that $(A,B)\in\tau$ implies $(A\cup A', B\cup B')\in\tau$: then also $(A',B')\in\tau$, since otherwise $(B',A')\in\tau$ and therefore $(B\cup B', A\cup A')\in\tau$.\looseness=-1

As $(A,B\cup B')\le (A,B)$, we have $(A,B\cup B')\in\tau$ by the consistency of~$\tau$. But $\{(A,B\cup B'), (B\cup B', A\cup A')\}\in\T_2$.%
   \COMMENT{}
   Therefore $(B\cup B', A\cup A')\notin\tau$, and hence $(A\cup A', B\cup B')\in\tau$ as desired.
  \endproof

\proclaim Lemma \xxxInfiniteB.
For every \AT\ $\tau$ and $(A,B)\in\tau$, the set $B$ is infinite.

\proof
If $B$ is finite, then $A\sim V$. Since $(B,V)\in\tau$ by Lemma~\xxxSmall, this implies $(B,A)\in\tau$ by Lemma~\xxxSymDiff.
   \endproof

\beginsection \secUFs. Tangles and ultrafilters

We are considering a fixed infinite graph $G=(V,E)$, with $\X$ denoting the set of finite subsets of~$V$. For each $X\in\X$, write $\C_X$ for the set of components of $G-X$, and let $\U_X$ denote the set of all \uf s on~$\C_X$.

Our aim in this section is to study those \AT s of $G$ that are not end tangles. We shall see that they define ultrafilters on the sets $\C_X$, and conversely that every way of choosing such ultrafilters consistently defines an \AT. More precisely, we shall see that the \AT s of~$G$ correspond to the points of a natural inverse limit of the sets~$\U_X$, with the end tangles among them corresponding to the limits of principal \uf s.

Recall that every end tangle $\tau = \tau_\omega$ of~$G$ defines a direction~$f$ in~$G$: a way of choosing for every $X\in\X$ one component of~$G-X$, the component $C$ in which $\omega$ lives. As we saw in Example~\xxxStarEx, an arbitrary \AT~$\tau$ may not select one $C\in\C_X$ in this way, but it still cannot orient the separations $\{A,B\}$ with $A\cap B = X$ arbitrarily. Indeed, as $\tau$ is consistent and has no subset in~$\T_{<\aleph_0}$,%
   \COMMENT{}
   it defines an ultrafilter on~$\C_X$, the collection of all $\C\sub\C_X$ that such that $\bigcup\C = B\sm A$ for some $(A,B)\in\tau$ with $A\cap B = X$:
 $$\textstyle U(\tau,X) := \big\{\,\C\sub\C_X\mid\big(\bigcup(\C_X\sm\C)\cup X, \bigcup\C\cup X\big) \in\tau\,\big\} \in\U_X.$$
 It is easy to check that this is indeed an ultrafilter. Every $(A,B)\in\tau$ with $A\cap B=X$ partitions $\C_X$ into two sets, the components in~$A$ versus those in~$B$, and $\tau$ puts the latter set in the \uf. The intersection of two such subsets of~$\C_X$ are also chosen by~$\tau$: if $(A,B), (A',B')\in\tau$ then also $(A\cup A', B\cap B')\in\tau$ by Lemma~\xxxP. By Lemma~\xxxSmall, these filter sets are non-empty, and they are closed under taking supersets in~$\C_X$ because $\tau$ is consistent.

For every $X\in\X$, we thus have a map
 $$\tau\mapsto U(\tau,X)\eqno\eqtauUtauX$$
 from $\Theta$ to~$\U_X$. These maps are not in general injective: if some $C\in\C_X$ is home to more than one end, for example, then $\{C\}$ will generate $U(\tau,X)$ for all the corresponding end tangles~$\tau$.%
   \COMMENT{}
   However, we shall see in Lemma~\xxxX\ that distinct $\tau,\tau'\in\Theta$ can never map to the same non-principal \uf\ as in~\eqtauUtauX.

We shall also see later in Lemma~\xxxSurjective\ that the maps in~\eqtauUtauX\ are nearly surjective: only principal \uf s in~$\U_X$ generated by $\{C\}$ for a {\it finite\/} component $C$ of~$G-X$ are not of the form~$U(\tau,X)$.%
   \COMMENT{}

The \uf s $U(\tau,X)$ for a given~$\tau$ but variable~$X$ are compatible for $X\sub X'$ just as the choices of $f(X)\in\C_X$ and $f(X')\in\C_{X'}$ are compatible when $f$ is a direction in~$G$. Indeed, consider the maps
 $$\fXX: \U_{X'}\to\U_X$$
 defined for all $X\sub X'\in\X$ by mapping an \uf\ $U'\!\in\U_{X'}$ to $U'\!\restricts X\sub 2^{\C_X}$, where $U'\!\restricts X$ is the set of all supersets%
   \COMMENT{}
   in~$\C_X$ of sets of the form~$\C'\!\restricts X$ with $\C'\!\in U'$ and $\C'\!\restricts X:= \{\,C\in\C_X\mid \exists C'\in\C': C'\sub C\,\}$. Less formally, from an \uf\ $U'\!\in\U_{X'}$ we obtain an \uf\ $U = \fXX(U')\in\U_X$ by putting the following sets $\C\sub\C_X$ in~$U$: pick some $\C'\!\in U'$, and let $\C$ consist of at least all those components of $G-X$ that contain some $C'\!\in\C'$ as a subset (Fig.~\figUFs).

 \figure \figUFs.
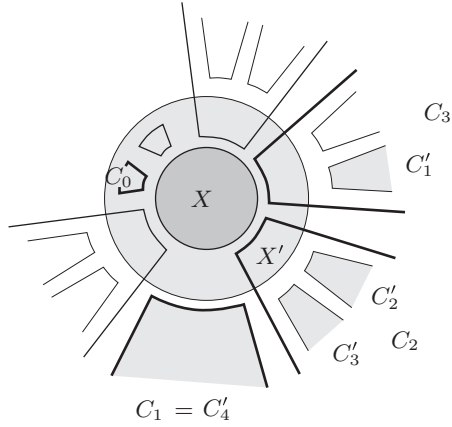
 $\!\C = \{C_0,C_1, C_2, C_3\} \in U$ if $\C'\! = \{C'_1,\dots,C'_4\}\in U'$
 (UFs; 1000)

\proclaim Lemma~\xxxInducedUFs.
The maps~$\fXX$ make $\,\UX$ into an inverse system, with $\X$ partially ordered by inclusion.

\proof
It is immediate that for $X\sub X'\sub X''$ we have $\fXX\circ f_{X''\!,X'} = f_{X''\!,X}$, as required for an inverse system, as long as these compositions are defined. But to establish this we first have to show that, given $X\sub X'$ and $U'\!\in\U_{X'}$, the set $U\!:= \fXX(U') \sub 2^{\C_X}$ is indeed an \uf.

To this end, notice that every component $C'$ of $G-X'$ lies inside some component $C$ of~$G-X$, because $C'$ is connected and does not meet $X\sub X'$. Let $f\: C'\mapsto C$ be this map from $\C_{X'}$ to~$\C_X$. Now consider a partition $\C_X = \C_1\cup\C_2$. For $i=1,2$ define $\C'_i := \{\,C'\!\in\C_{X'}\mid f(C')\in\C_i\,\}$. Since $f$ has domain all of~$\C_{X'}$, the $\C'_i$ partition~$\C_{X'}$. As $U'$ is an \uf\ on~$\C_{X'}$, it contains exactly one of the~$\C'_i$ as an element. By definition of~$\fXX$, the corresponding~$\C_i$%
   \COMMENT{}
   lies in~$U$: it is a superset%
  \Footnote{The superset can be strict if $G-X$ has (finite) components contained in~$X'$.}
   of~$\C'_i\restricts X = \{\,f(C')\mid C'\in\C'_i\,\}$.

The other two properties required of an \uf, that it is closed under taking finite intersections and does not contain~$\es$, follow for~$U$ from the corresponding properties of~$U'$ in a similar way.%
   \COMMENT{}
   \endproof

For want of a better expression, let us call the ultrafilters in the sets~$\U_X$ the {\it \uf s of cofinite components\/} in~$G$, and put
 $$\U := \varprojlim (\,\U_X\mid X\in\X\,).$$

Conveniently, our maps $\tau\mapsto U(\tau,X)$ commute with the maps~$\fXX$:

\proclaim Lemma \xxxCommute.
For all $X\sub X'\in\X$ and all \AT s~$\tau$ we have
 $$\fXX \big(U(\tau,X')\big) = U(\tau,X).$$
 Every \AT~$\tau$ therefore defines a limit $\upsilon_\tau = (\,U_X\mid X\in\X\,)\in\,\U$ in which $U_X = U(\tau,X)$ for all~$X$.

\proof
Let $U\!:= \fXX\big( U(\tau,X')\big)$, consider any $\C\in U$, and let $(A,B)\in\vS$ be such that $A\cap B = X$ and $B\sm A = \bigcup\C$. Our aim is to show that $(A,B)\in\tau\,$: then $\C\in U(\tau,X)$, giving $U\sub U(\tau,X)$, with equality since both are \uf s.

By definition of~$\fXX$, there exists $\C'\!\in U(\tau,X')$ such that $\C\supe \C'\!\restricts X$. The separation $(A',B')\in\vS$ with $A'\cap B' = X'$ and $B'\sm A' = \bigcup\C'$ thus lies in~$\tau$. Since $(A,B\cup X')\le (A',B')$,%
   \COMMENT{}
   this implies $(A,B\cup X')\in\tau$ by the consistency of~$\tau$. Now $(A,B)\in\tau$ follows by Lemma~\xxxSymDiff.
 \endproof

Conversely, every $\upsilon\in\U$ comes from a tangle in this way:

\proclaim Lemma \xxxAllTangles.
For every limit $(\,U_X\mid X\in\X\,)\in\,\U$ there exists a unique \AT\ $\tau$ in~$G$ such that $U_X = U(\tau,X)$ for all~$X\in\X$. The map
 $$\tau\mapsto\upsilon_\tau\eqno\eqtauy$$
 defined in Lemma \xxxCommute, therefore, is a bijection from $\Theta$ to $\U$.

\proof
Let $\tau := \{\,(A,B)\in\vS\mid \C_{A\cap B}\cap 2^B \in U_{\!A\cap B}\,\}$. Clearly, $U(\tau,U)=U_X$ for all $X\in\X$ if $\tau$ is indeed an \AT, so let us check this; uniqueness will be clear, since distinct \AT s $\tau,\tau'$ differ on some separation $\{A,B\}$, so that $U(\tau,X)\ne U(\tau', X)$ for $X=A\cap B$.%
   \COMMENT{}

$\!$For every separation $\{A,B\}\in S$, the sets $\C_X\cap 2^A$ and $\C_X\cap 2^B$ partition~$\C_X$, so $U_X$ contains exactly one of them. Hence, $\tau$ is an orientation of~$S$.

To show that $\tau$ is consistent, consider $(A,B) < (A',B')\in\tau$. Let $X:= A\cap B$ and $X':= A'\cap B'$, and put $X'':= X\cup X'$. Let $\C:= \C_X\cap 2^A$ and $\C':=\C_{X'}\cap 2^{B'}\!$. Note that $\C,\C'\sub \C_{X''}$, since $(A,B)\le (A',B')$. Now if $(B,A)\in\tau$, then
 $$f_{X'',X}\! : U_{X''}\mapsto U_X\owns\C\supe\C''\restricts X$$
 for some $\C''\in U_{X''}$. Since $\C\sub\C_{X''}$, this means that in fact $\C\supe\C''$,%
   \COMMENT{}
   and hence that $\C\in U_{X''}$. Similarly, $\C'\in U_{X''}$. But $\C\cap\C'=\es$, a contradiction.

\goodbreak

It remains to show that $\tau$ has no subset $\sigma=\{(A_1, B_1), (A_2, B_2), (A_3, B_3)\}$ in~$\T_3$. Suppose it does. For $i=1,2,3$ let $X_i:= A_i\cap B_i$, and put $X:= \bigcap_i B_i$. This is finite by definition of~$\T_3$, and includes every $X_i$ since $\sigma$ is a star. Then
 $$f_{X,X_i}\!: U_X\mapsto U_{X_i}\owns \C_{X_i}\cap 2^{B_i}\supe\C_i\restricts X_i$$
 for some $\C_i\in U_X$, for each~$i$. Every $C\in\C_i$ is a subset of some $C_i\in \C_{X_i}\cap 2^{B_i}$, and hence of~$B_i$. Hence any $C\in\C_1\cap\C_2\cap\C_3$ is both a subset of $\bigcap_i B_i = X$ and an element of~$\C_X$, which is impossible. As $\C_1\cap\C_2\cap\C_3\in U_X$ is non-empty, this is a contradiction.
 \endproof

We have seen that every end $\omega$ of~$G$ defines a tangle~$\tau_\omega$ by~\eqomegatau, which in turn defines a limit $\upsilon\in\U$ by~\eqtauy. The composition
 $$\omega\mapsto\tau_\omega\mapsto\upsilon_{\tau_\omega}\eqno\eqomegay$$
 maps $\omega$ to a limit $\upsilon = (\,U_X\mid X\in\X\,)$ in which every $U_X$ is a principal \uf\ in~$\U_X$, generated by~$\{C\}$ say, where $C$ is the component of $G-X$ in which $\omega$~lives. The converse of this is also true: if every $U_X$ in $(\,U_X\mid X\in\X\,) = \upsilon_\tau$ is principal, then $\tau$ is an end tangle:

\proclaim Lemma~\xxxPrincipal.
The following assertions are equivalent for all \AT s~$\tau$:
 \smallskip
 \plainitem{$\bullet$} $\tau$ is an end tangle;
 \smallskip
 \plainitem{$\bullet$} $U(\tau,X)$ is a principal \uf, for all $X\in\X$.
 \enditem\medskip\medskip

\proof
We only have to show the backward implication. For each $X\in\X$, let $C_X$ be the unique component of $G-X$ such that $\{C_X\}$ generates~$U(\tau,X)$. By Lemma~\xxxCommute, the map $f\: X\mapsto C_X$ is a direction in~$G$. Since the map $\omega\mapsto f$ is surjective~[\the\refTopEnds], there is an end~$\omega$ of~$G$ such that $f = f_\omega$. This $\omega$ lives in every~$C_X$, giving $\tau = \tau_\omega$ as desired.
 \endproof

Let us call $\tau\in\Theta$ an {\it ultrafilter tangle\/} if at least one of the \uf s~$U_X$ in $\upsilon_\tau = (\,U_X\mid X\in\X\,)$ is non-principal. Of every such $X$ we say that it {\it witnesses\/} that $\tau$ is an \uf\ tangle.

Lemma~\xxxPrincipal\ tells us that the \AT s in $G$ divide into its end tangles and its ultrafilter tangles. We have thus proved Theorem~\xxxUFs:

\proclaim Theorem~\xxxUFs. Let $G$ be any graph.
\pitem{i} The \AT s of $G$ are precisely the limits of the inverse system of its sets of \uf s of cofinite components. 
\pitem{ii} The ends of $G$ are precisely those of its \AT s whose \uf s of cofinite components are all principal.
 \enditem

\beginsection \secUFTs. A closer look at ultrafilter tangles

Recall that our aim was to understand better those \AT s that are not defined by an end. We have seen that these are the ultrafilter tangles, those $\tau\in\Theta$ such that at least one of the $U_X = U(\tau,X)$ in their $\upsilon_\tau = (\,U_X\mid X\in\X\,)$ is non-principal. In this section we show that each of these $U_X$ already determines~$\tau$. Thus, every ultrafilter tangle is determined by a single \uf\ of co\-finite components, not just by a limit of such \uf s.

In particular, if $U = U(\tau,X)$ is non-principal and $X\sub X'$, then both $\upsilon_\tau$ and $U' = U(\tau,X')$ are determined by~$U$, since $U$ determines~$\tau$ and $\tau$ determines $\upsilon_\tau$ and~$U'$. However we shall prove this directly, without involving~$\tau$: we show that for every element $U$ of the set
 $$\U^*_X := \{\,U\in\U_X\mid U\hbox{ is non-principal}\,\}$$
 there is a unique $U'\!\in\U_{X'}$ such that $\fXX(U')=U$. (This $U'$ will also lie in~$\U^*_{X'}$.) Thus, the maps $\fXX$ have inverses on the sets $\U^*_X$ of non-principal \uf s. Hence%
   \COMMENT{}
   there is also a unique $\upsilon\in\U = \varprojlim (\,U_X\mid X\in\X\,)$ with $U_X = U$.%
   \COMMENT{}

We finally show that the set $\X_\tau$ of all $X\!\in\X$ witnessing that a given~$\tau$ is an ultrafilter tangle has a least element~$X$. From its corresponding $U = U_X$ we can thus directly construct all the \uf s of cofinite components induced by~$\tau$, the filters $U(\tau,X')$ with $X'\!\in\X_\tau$, by applying the inverses of the maps~$\fXX$.%
   \COMMENT{}

\proclaim Lemma \xxxg.
Let $X\sub X'\in\X$, and let $U\!\in\U^*_X$ be a non-principal \uf\ on~$\C_X$. Then there is a unique \uf\ $U'\!$ on~$\C_{X'}$ such that $\fXX(U') = U$. This \uf~$U'\!\in\U_{X'}$ is also non-principal, and it satisfies
 $$U' = \{\,\C\sub\C_{X'}\mid \exists \D\in U\:\ \D\sub\C\,\}.$$
 \smallskip

\proof
As $\fXX$ maps principal \uf s to principal \uf s, it is clear that any $U'\!\in\U_{X'}$ satisfying $\fXX(U')=U$ lies in~$\U^*_{X'}$. Let us show that there is a unique such~$U'\!\in\U_{X'}$.

Consider any bipartition $\C_{X'} = \C\cup\C'$. Since every component of $G-X$ that does not meet~$X'$ is also a component of~$G-X'$, our partition of~$\C_{X'}$ induces a partition $\D^-\cup\D\cup\D'$ of~$\C_X$, where $\D^-$ is the set of components of $G-X$ meeting~$X'$, and $\D\sub\C$ and $\D'\sub\C'$.

As $U$ is non-principal, it does not contain the finite set~$\D^-$. Hence exactly one of $\D$ and~$\D'$ lies in~$U$, say~$\D$.%
   \COMMENT{}
   Then any $U'\!\in\U_{X'}$ satisfying $\fXX(U')=U$ contains~$\C$: otherwise it would contain~$\C'$, and then $U_X$ would contain $\D^-\cup\D'\supe\C'\restricts X$ which it does not.%
   \COMMENT{}

Let $U'$ be the set of all $\C$ obtained in this way: the set of all $\C\sub\C_{X'}$ such that $\C\cap\C_X\!\in U$. If $U'$ is a filter, it will be an ultrafilter, because it contains a set from every bipartition of~$\C_{X'}$. It thus remains to show that $U'$ is indeed a filter, and that it satisfies $\fXX(U') = U$; we have already seen that $U'$ will then be unique in~$\U_{X'}$ with this property.%
   \COMMENT{}

 Since every $\C\in U'$ has a subset $\D\in U$, clearly $\es\notin U'$. And for $\C_1,\C_2\in U'$, with $\D_i\sub\C_i$ in~$U$ say, we have $U\!\owns\D_1\cap\D_2\sub\C_1\cap\C_2\cap\C_X$, giving $\C_1\cap\C_2\in U'$.%
   \COMMENT{}%
   \COMMENT{}
   Thus, $U'\in\U_{X'}$.

It is straightforward from the definitions that $\fXX(U')\sub U$.%
   \COMMENT{}
   Since both these are \uf s on the same set (cf.\ Lemma~\xxxInducedUFs), this implies the desired equality.
   \endproof

As a consequence of Lemma~\xxxg, we have inverses of the maps $\fXX$ on the non-principal \uf s:

\proclaim Lemma \xxxgf.
For all $X\sub X'\in\X$ there exists a map $g_{X,X'}\:\> \U^*_X\to\U_{X'}$ such that $\fXX\circ g_{X,X'}$ is the identity on~$\U^*_X$.\noproof

Our aim was to show that if $\tau\in\Theta$ is an \uf\ tangle and $X$ witnesses this, then $U(\tau,X)$ alone determines $\tau$. This is an easy consequence of Lemma~\xxxg:

\proclaim Lemma~\xxxX.
Let $\tau\in\Theta$ and $X\in\X$ be such that $U(\tau,X)$ is non-principal. Then $U(\tau',X)\ne U(\tau,X)$ for every $\tau'\in\Theta\sm\{\tau\}$.

\proof
As $\tau\ne\tau'$, there exists $(A,B)\in\tau$ such that $(B,A)\in\tau'$. For $X':= A\cap B$ this gives $U(\tau,X')\ne U(\tau',X')$. Then also $U(\tau, X\cup X')\ne U(\tau',X\cup X')$: if these were the same filter $U\!\in\U_{X\cup X'}$ we would have $U(\tau,X') = f_{X\cup X'\!,X'}(U) = U(\tau',X')$ by Lemma~\xxxCommute. By Lemma~\xxxg, the fact that $U(\tau, X\cup X')\ne U(\tau',X\cup X')$ implies $U(\tau',X)\ne U(\tau,X)$.%
   \COMMENT{}
   \endproof

Every $U\!\in\U^*_X$ can be used to define the unique $\tau$ with $U(\tau,X) = U$ directly, by telling us which orientation of a given separation $\{A,B\}\in S$ lies in~$\tau$:

\proclaim Lemma \xxxDirect.
Let $\tau\in\Theta$ be an \uf\ tangle, witnessed by $X\in\X$. Then
 $$\textstyle\tau = \big\{\, (A',B')\in\vS\mid \exists\C\in U(\tau,X): \bigcup\C\sub B'\,\big\}.$$
 \smallskip

\proof
For a proof of~`$\sub$'%
   \COMMENT{}
   let $(A',B')\in\tau$ be given, and put $X' := A'\cap B'$. Pick $(A,B)\in\tau$ with $A\cap B = X$. By Lemma~\xxxP, also $(A'',B'')\in\tau$ for $A''=A\cup A'$ and $B'' = B\cap B'$.%
   \COMMENT{}
   For $X'' = X\cup X'$ then also $(A'',X''\cup B'')\le (A'',B'')$ lies in~$\tau$, by the consistency of~$\tau$. The set $\C''$ of components of $G-X''$ contained in~$B''$ then lies in~$U(\tau, X'')$. As $f_{X'',X}(U(\tau,X'')) = U(\tau,X)$ by Lemma~\xxxCommute, the set $\C''$ has a subset $\C\in U(\tau,X)$ by Lemma~\xxxg. Then $\bigcup\C\sub\bigcup\C''\sub B''\sub B'$ as desired.

Conversely, any $(A',B')\in\vS$ with $\bigcup\C\sub B'$ for some $\C\in U(\tau,X)$%
   \COMMENT{}
   must lie in~$\tau$: otherwise $(B',A')\in\tau$ with $\bigcup\C'\sub A'$ for some $\C'\in U(\tau,X)$, as shown above, but $\C\cap\C'$ is finite and hence not in~$U(\tau,X)$, a contradiction.
   \COMMENT{}
 \endproof

Note that Lemma~\xxxDirect\ also implies Lemma~\xxxX. Also, we could replace the requirement of $\bigcup\C\sub B$ with $\bigcup\C\sub B\sm A$: our proof yield this directly, but it also follows retrospectively since $U$ is non-principal and $A\cap B$ can meet only finitely many elements of~$\C$.

\medbreak

For an end tangle $\tau = \tau_\omega$, a single ray $R\in\omega$ can be used as an `oracle' to determine which of the two orientations of a given separation $\{A,B\}\in S$ is in~$\tau$: it is the unique orientation that points to a tail of~$R$. Lemma~\xxxDirect\ says that for \uf\ tangles we have similar oracles, given by a single \uf\ in
 $$\U_\X:= \bigcup_{X\in\X}\U_X.$$

We have thus proved the following more precise version of Theorem~\xxxTau:

\proclaim Theorem \xxxtau.
Every \AT\ $\tau$ in~$G$ satisfies exactly one of the following:
 \medskip
 \plainitem{$\bullet$} $\exists$ ray $R\sub G$ such that
    $\tau = \{\,(A,B)\in\vS : G[B]\hbox{ contains a tail of }R\,\}${\rm;}
 \smallskip
 \plainitem{$\bullet$} $\exists$ \uf\ $U\!\in\U_\X$ such that
    $\tau = \{\,(A,B)\in\vS\mid \exists\C\in U\,\:\>\bigcup\C\sub B\sm A\,\}$.%
   \COMMENT{}
   \smallskip\medskip

\proof
We have shown everything claimed, except that no tangle can satisfy both statements. But this is clear: an \uf\ $U$ as in the second statement cannot be principal since, given its generating set~$\{C\}$, we can choose $\{A,B\}\in S$ with $A\cap B\cap C\ne\es$, in which case $\bigcup\C\supe C$ will not be contained in either $B\sm A$ or~$A\sm B$, so $\tau\cap \{(A,B),(B,A)\} = \es$, a contradiction.
 \endproof

Let us now prove that, for every ultrafilter tangle $\tau\in\Theta$, the set
 $$\X_\tau:= \{\,X\in\X\mid U(\tau,X)\in\U^*_X\,\}$$ 
 of all $X$ witnessing that $\tau$ is an ultrafilter tangle is the up-closure in $\X$ of a single element:

\proclaim Theorem~\xxxXtau.
For every \uf\ tangle~$\tau$ the set $\X_\tau$ has a least element~$X_\tau$. Then $\X_\tau = \{\,X\in\X\mid X_\tau\sub X\,\}$.

\proof
We already noted in Lemma~\xxxg\ that $\X_\tau$ is closed upwards in~$\X$. It remains to show that it has a least element.

Suppose not. Let $X',X''$ be incomparable minimal elements of~$\X_\tau$. Pick $x'\in X'\sm X''$,%
   \COMMENT{}
   and let $X:= X'\sm\{x'\}$. By the minimality of~$X'$, the \uf\ $U(\tau,X)$ is principal, and hence generated by $\{C\}$ for some $C\in\C_X$.%
   \COMMENT{}
   As $U(\tau,X')$ is non-principal and $\fXX(U(\tau,X')) = U(\tau,X)$, the set $\C$ of components of $C-x'$ lies in~$U(\tau,X')$.%
   \COMMENT{}
   As $X''$ meets only finitely many elements of~$\C$, the others form a set $\C'\!\in U(\tau,X')$.%
   \COMMENT{}
   Similarly, pick $x''\!\in X''\sm X'$ and find a set $\C''\!\in U(\tau,X'')$ of components of $G-X''$ that avoid~$X'$.

Every $C'\in\C'$ lies inside the same component $C''$ of $G-X''$ as~$x'$, because both avoid~$X''$ but $G$ contains an edge from $x'$ to~$C'$: otherwise, $C'$~would be in~$\C_X$, which it is not since $C-x'$ contains it.%
  \COMMENT{}
   Since the components in $\C''$ avoid $X'\!\owns x'$, we thus have $C'\sub C''\!\notin\C''$ for every $C'\in\C'$.

Thus, $\bigcup\C'$ and $\bigcup\C''$ are disjoint sets of vertices separated by the finite set~$X''$ (as well as by~$X'$);%
   \COMMENT{}
   let $\{A,B\}\in S$ with $A\cap B = X''$ separate them. Both orientations of $\{A,B\}$ must be in~$\tau$, by Lemma~\xxxDirect\ applied with $X'$ and with~$X''$, respectively. This contradicts our assumption that $\tau$ is an orientation of~$S$.
   \endproof

Finally, let us go back and use the maps $g_{X,X'}$ from Lemma~\xxxgf\ to prove that the maps $\tau\mapsto U(\tau,X)$ in~\eqtauUtauX\ are essentially surjective; we shall need this in our proof of Theorem~\xxxTop.

\proclaim Lemma \xxxSurjective.
Let $X\in\X$, and let $U\!\in\U_X$ be an \uf\ on~$\C_X$ not generated by $\{C\}$ for any finite component $C$ of $G-X$. Then there exists an \AT\ $\tau\in\Theta$ such that $U = U(\tau,X)$.

\proof
Assume first that $U$ is non-principal. Our aim is to find a limit point $\upsilon = (\,U_Y\mid Y\in\X\,)\in\U$ such that $U_X = U$. By Lemma~\xxxAllTangles\ there will then exist some $\tau\in\Theta$ with $\upsilon = \upsilon_\tau$, for which $U(\tau,X) = U_X = U$ as desired.

 \figure \figBigDiagram.
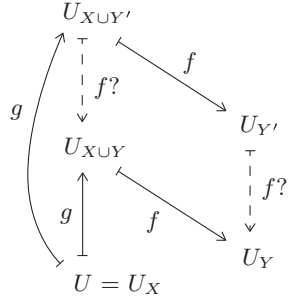
 The known maps are drawn as solid lines, the desired maps as broken lines
 (BigDiagram; 1000)

For every $X'\in\X$ with $X\sub X'$ let $U_{X'} := g_{X,X'} (U)$, and for all other $Y\!\in\X$ let $U_Y:= f_{X',Y} (U_{X'})$ for $X':= X\cup Y$ (Fig.~\figBigDiagram). Then, in fact,
 $$U_Y = f_{X\cup Y,Y} \circ g_{X,X\cup Y} (U)$$
 for every $Y\in\X$, and $U_X = U$. To show that $(\,U_Y\mid Y\in\X\,)\in \U$, we have to show that $f_{Y',Y} (U_{Y'}) = U_Y$ for all $Y\sub Y'\in\X$. To see this,  note first that
 $$\eqalign{U_{X\cup Y} &= g_{X,X\cup Y}(U)\cr
       &= g_{X,X\cup Y}(f_{X\cup Y', X}(U_{X\cup Y'}))\cr \COMMENT{}
       &= g_{X,X\cup Y}(f_{X\cup Y, X}(f_{X\cup Y',X\cup Y}(U_{X\cup Y'})))\cr
       &= f_{X\cup Y',X\cup Y}(U_{X\cup Y'}).\cr}$$ \COMMENT{}

\noindent
   Hence%
   \COMMENT{}
 $$\eqalign{f_{Y',Y}(U_{Y'}) &= f_{Y',Y} (f_{X\cup Y',Y'} (g_{X,X\cup Y'} (U)))\cr
   &= f_{Y',Y} (f_{X\cup Y', Y'} (U_{X\cup Y'}) )\cr
   &= f_{X\cup Y', Y} (U_{X\cup Y'})\cr
   &= f_{X\cup Y,Y} (f_{X\cup Y'\!,\, X\cup Y} (U_{X\cup Y'}) )\cr
   &= f_{X\cup Y,Y} (U_{X\cup Y})\cr
   &= f_{X\cup Y,Y} (g_{X,X\cup Y}(U))\cr
   &= U_Y,\cr}$$
 as desired.

Suppose now that $U$ is principal, generated by $\{C\}$ with $C\in\C_X$ say, and that $C$ is infinite. If $C$ contains a ray, the end~$\omega$ of this ray defines $\tau_\omega$, for which $U = U(\tau_\omega,X)$ as desired. If not, then $C$ has a finite set~$Z$ of vertices such that $C-Z$ has infinitely many components: otherwise we could construct a ray $z_0 z_1\dots$ in $C$ inductively by choosing each $z_n$ from an infinite component of $C-\{z_0,\dots,z_{n-1}\}$.%
   \COMMENT{}
   Pick a non-principal \uf~$U_Z$ on the set of components of~$C-Z$, notice that these are also components of~$G-X'$ for $X' = X\cup Z$,%
   \COMMENT{}
   and let $U'$ be the (non-principal) \uf\ on~$\C_{X'}$ generated by~$U_Z$.%
   \COMMENT{}
   Then $\fXX(U')=U$, by definition of $U'$ and~$\fXX$.%
   \COMMENT{}
   By the case already treated, there exists $\tau\in\Theta$ such that $U' = U(\tau, X')$. This $\tau$ achieves our aim,%
   \COMMENT{}
   since 
 $$U(\tau,X) = \fXX(U(\tau,X')) = \fXX(U') = U$$
 by Lemma~\xxxCommute.%
   \COMMENT{}
 \endproof

\beginsection \secTop. Tangles at infinity: compactifying an arbitrary graph

Our aim in this section is to prove Theorem~\xxxTop: that the \AT s can be used as points at infinity to compactify~$G$, in a way that yields its Freudenthal compactification when $G$ is locally finite (and hence all its \AT s are end tangles). To make this process more transparent we shall first define a topology on~$\U$ itself, which can be done in a rather canonical way. We shall then adapt this to define a topology on all of $G\cup\U$ that induces this topology on~$\U$, as well as the usual 1-complex topology on~$G$. This space $|G| = G\cup\U$ or, equivalently by Lemma \xxxAllTangles, the space $|G| = G\cup\Theta$, will be the desired compactification of~$G$.%
   \looseness=-1

For each $X\in\X$, take the topology on~$\U_X$ whose basic open sets are those of the form%
   \COMMENT{}
 $$\O(\C) := \{\,U\!\in\U_X\mid \C\in U\,\},$$
one for each $\C\sub\C_X$. This topology, the Stone topology on~$\C_X$, makes $\U_X$ into a compact topological space%
   \Footnote{The \SC\ compactification of~$\C_X$. Its compactness is immediate from Tychonov's theorem when we view $\U_X$ as a subspace of~$2^{2^{\C_X}}\!$: the set $\U_X$ is closed in this compact space, since any violation of the \uf\ axioms involves only finitely many subsets of~$\C_X$.}%
   \COMMENT{}
   in which the principal \uf s on $\C_X$ are dense.%
   \COMMENT{}%
   \penalty-200\
   The space $\U_X$ is clearly Hausdorff, indeed totally disconnected.%
   \COMMENT{}

Our inverse system $\UX$ is compatible with these topologies:%
   \COMMENT{}

\proclaim Lemma \xxxcts.
The maps $\fXX\:\>\U_{X'}\to\U_X$ are continuous.
 \noproof%
   \COMMENT{}

\noindent
In fact, for $X\sub X'\in\X$ the open sets $\fXX^{-1}(\O(\C))\sub\U_{X'}$ are themselves basic:
 $$\fXX^{-1}(\O(\C)) = \O(\C')\quad{\rm for}\quad\C' = \{\,C'\in\C_{X'}\mid \exists C\in\C\,\:\> C\supe C'\,\}.\eqno\eqbasic$$%
   \COMMENT{}

Topologizing the $\U_X$ has an interesting windfall for graphs without ends:%
   \COMMENT{}

\proclaim Proposition \xxxTanglesExist.
Every infinite graph has an \AT.

\proof
Inverse limits of non-empty compact spaces are non-empty~[\the\refInvLim]. Hence $\U$ is not empty, and so by Lemma~\xxxAllTangles\ neither is~$\Theta$.
   \endproof

Let us give the set $\Theta=\U = \varprojlim\UX$ of \AT s in~$G$ the subspace topology from $\prod_{X\in\X}\U_X$ endowed with the product topology of the~$\U_X$.%
   \COMMENT{}
   \looseness=-1

\proclaim Proposition \xxxTheta.
The topological space of all \AT s of an infinite graph is compact and totally disconnected.%
   \COMMENT{}

\proof
Since the $\U_X$ are Hausdorff and the $\fXX$ are continuous, the space $\U = \varprojlim (\,\U_X\mid X\in\X\,)$ is closed in~$\prod_{X\in\X}\U_X$,%
   \COMMENT{}
   which inherits its own compactness from that of the~$\U_X$ by Tychonov's theorem. It is totally disconnected, because the~$\U_X$ are.%
   \COMMENT{}
   \endproof

To describe this topology more explicitly, consider any $X\in\X$ and $\C\sub\C_X$. Let
 $$\O(X,\C) := \{\,\upsilon\in\U\mid \C\in U_X\,\} = f^{-1}_X (\O(\C)),$$
 where $\upsilon = (\,U_X\mid X\in\X\,)$ and $f_X$ is the (continuous) projection $\U\to\U_X$.%
   \COMMENT{}

\proclaim Lemma \xxxBasis.
The sets $\O(X,\C)$%
   \COMMENT{}
   form a basis of open sets in~$\U$.

\proof
   \COMMENT{}
   By definition of the topology on~$\U$, these sets form a subbasis. By~\eqbasic\ they  even form a basis, because every finite intersection of such sets $\O(X,\C)$ can be rewritten as the union of sets $\O(X',\C')$ with $X'$ the (finite) union of these~$X$.%
   \COMMENT{}
   \endproof%
   \COMMENT{}

When $G$ is locally finite and connected, it is compactified by its ends in the so-called {\it Freudenthal compactification\/}~[\the\refFreudenthal,\th\the\refFreudenthalNeu;\th\the\refBook]. The following definition for our arbitrary~$G$ defaults to this when $G$ is locally finite and connected,%
   \COMMENT{}
   and hence all its \AT s are end tangles.

Let us view $G$ as a 1-complex with the usual topology. Its edges are copies of the real interval~$[0,1]$, and choosing for every edge $e=vw$ at some vertex~$v$ any%
   \COMMENT{}
   half-open partial edge $e_x=[v,x)\sub e$ with $x\in\interior e$, makes $\bigcup_e e_x$ into an open neighbourhood of~$v$. Let us extend $G$ to a topological space
 $$|G| = G\cup\, \U = G\cup\Theta$$
(cf.\ Lemma \xxxAllTangles) by also%
   \COMMENT{}
   declaring as open, for all $X\!\in\X$ and all $\C\sub\C_X$,%
   \COMMENT{}
   the sets
 $$\textstyle\O_G (X,\C) := \bigcup\C \,\cup\ \interior E(X,\bigcup\C)\,\cup\  \O(X,\C)\eqno
\eqopensets$$
 and taking the topology on $|G|$ that this generates.%
   \COMMENT{}
   Here, $\interior E(X,\bigcup\C)$ is the set of all inner points of edges between $X$ and components of $G-X$ in~$\C$. Note that the subspace topology on $\U\sub|G|$ is our original topology on~$\U$, and that the subspace topology on~$G$ is its orginal 1-complex topology.%
   \COMMENT{}

Let us prove that $|G| = G\cup\U = G\cup\Theta$ is a compact space, the {\it tangle compactification\/} of~$G$:

\proclaim Theorem~\xxxTop. Let $G$ be any graph.
\pitem{i} $|G|$ is a compact space in which $G$ is dense and $|G|\sm G$ is totally disconnected.
\pitem{ii} If $G$ is locally finite and connected, then all its \AT s are ends, and $|G|$ coincides with the Freudenthal compactification of~$G$.
   \enditem\medskip\medskip

\proof (i) Consider any cover $\O$ of $|G|$ by open subsets of~$G$ and basic open sets of the form~$\O_G(X,\C)$. Since $\U$ is compact in the subspace topology of~$|G|$, this has a finite subset of the form
 $$\F = \{\,\O_G (X,\D_X)\mid X\in\X'\,\}$$
(with $\X'\sub\X$ finite) that covers~$\U$. Our aim is%
   \COMMENT{}
   to show that, for $X'\!:=\bigcup\X'$, the sets in~$\F$ and $G[X']$ together cover~$|G|$: since $G[X']$ is a finite graph and hence compact, there will then also be a finite subcover of~$\O$ for all of~$|G|$.

For every $X\in\X'$, let
 $$\D'_X := \{\,C'\in\C_{X'}\mid \exists C\in\D_X\:\> C\supe C'\,\}.$$%
   \COMMENT{}
    Then $\fXX^{-1}(\O(\D_X)) = \O(\D'_X)$ by~\eqbasic. As $f_X = \fXX\circ f_{X'}$, this implies%
   \COMMENT{}
 $$\O(X,\D_X) = f_X^{-1}(\O(\D_X)) = f_{X'}^{-1}(\O(\D'_X)) = \O(X',\D'_X)\eqno\equps$$
 and hence $\O_G (X,\D_X)\supe \O_G (X',\D'_X)$, since $\bigcup\D_X\supe\bigcup\D'_X$ by definition of~$\D'_X$.%
   \COMMENT{}

\goodbreak

It thus suffices to show that the $\O_G (X',\D'_X)$ cover~$|G|\sm G[X']$, i.e., that every component of $G-X'$ is an element of some $\D'_X$ with $X\in\X'$. Suppose not, and let
 $$\C' :=\ \C_{X'}\sm\bigcup_{X\!\in\X'}\D'_X.$$%
   \COMMENT{}
If $\bigcup\C'$ is a finite graph,%
   \COMMENT{}
   we add this to $G[X']$ and achieve $\C'=\es$ as desired. We may therefore assume that $\bigcup\C'$ is infinite.

If $\C'$ contains an infinite component~$C'$, let $U'\!\in\U_{X'}$ be the \uf\ on $\C_{X'}$ generated by~$\{C'\}$. If not, then $\C'$ is infinite; pick a non-principal \uf\ on~$\C'$ and let $U'\!\in\U_{X'}$ be the (non-principal) \uf\ it generates on all of~$\C_{X'}$.%
   \COMMENT{}
   By Lemmas \xxxSurjective\ and \xxxAllTangles, there exists $\upsilon = (\,U_X\mid X\in\X\,)\!\in\U$ such that $U_{X'} = U'$. By~\equps\ and the definition of~$\C'$, this $\upsilon$ does not lie in~$\bigcup\F$, a contradiction.

 (ii) This is easy; see~[\the\refBook] for the definition of the Freudenthal compactification. To see that the topology for $|G|$ defined there coincides with ours here, remember that if $G$ is locally finite then any finite $X\sub V$ sends only finitely many edges to $G-X$.%
   \COMMENT{}
   \endproof

We remark that, as defined above, $|G|$ is not in general Hausdorff: the centre of an infinite star, for example, cannot be topologically separated from any open set containing a tangle of that star. However, every two points of $V\cup\Theta\sub |G|$ have neighbourhoods in $|G|$ that meet only in inner points of edges. If we delete the set $\interior E$ of all inner edge points from~$|G|$, the resulting space $|G|\sm\interior E$ will be a Hausdorff compactification of~$V$ that still reflects the structure of~$G$.%
   \COMMENT{}

Conversely, if we are prepared to give up the compactness of~$|G|$ (while keeping $\U = |G|\sm G$ and $V\cup\Theta = |G|\sm\interior E$ compact, which may be more crucial), we can make $|G|$ itself Hausdorff: we just have to allow more open sets in~\eqopensets\ by replacing $\interior E(X,\bigcup\C)$ with unions of either arbitrary half-edges~$(y,z)\sub (x,z)$ for each $e=(x,z)\in\interior E(X,\bigcup\C)$, or by uniformly chosen such half-edges (where all $y$ have distance some fixed positive~$\epsilon<1$ from $z\in\bigcup\C$ when the edge $e = [x,z]$ is viewed as a copy of the real interval~$[0,1]$.

\beginsection \secBlocks. Closed tangles

When Robertson and Seymour introduced tangles for finite graphs, their intended key feature was that they point to parts of the graph that are in some sense highly connected. Any large enough grid, for example, defines a $k$-tangle for any fixed~$k$, even though it is not highly connected as a subgraph: since every separation $\{A,B\}$ of order~$<k$ leaves most of the grid on one side, it can be oriented `towards' that side, and these orientations satisfy the tangle axioms.

\goodbreak

Our \AT s do not all point to a highly connected part of~$G$. Indeed, $G$~could be a locally finite tree, but it would still have end tangles pointing to its ends~-- which can hardly be seen as highly connected structures in any sense. On the other hand, an infinite complete subgraph also defines an \AT, for which a better case could be made. Ultrafilter tangles, however, have no connected~-- let alone highly connected~-- focus at all.

Some attempts have been made to at least identify those kinds of ends that tell us where our graph is highly connected. Candidates included the Halin ends that are not Freudental (or {\it topological\/}) ends~[\the\refCarmesinEF], which are those that have one or more vertices send an infinite fan to each of their rays~[\the\refTopEnds]. The earliest attempt, perhaps, was to consider `thick' ends~[\the\refHalinGrid,\th \the\refCTW], those of infinite (vertex) degree: these are the ends that contain an infinite set of disjoint rays~[\the\refBook,\th \the\refTopSurvey]. Halin~[\the\refHalinGrid] showed that these are precisely the ends (whose \AT\ is) defined by a half grid minor. An obvious analogue would be to consider the ends defined by a full grid minor~-- these have been characterized by Heuer~[\the\refHeuerFullGrid]~--%
   \COMMENT{}
   or infinite clique minors or subdivision as in~[\the\refRST,\th\the\refRSTminors].

I would like to propose a new alternative: that an \AT\ is deemed to signify a highly connected part of~$G$ if and only if it is closed in a certain natural topology on~$\vS = \vS_{\A}$.%
   \COMMENT{}
   We shall be able to characterize those tangles in graph-theoretical terms. They will all be end tangles, including those defined by an infinite complete subgraph but not, for example, the end tangles of a tree.\looseness=-1

The topology on~$\vS$ has the following basic open sets. Pick a finite set $Z\sub V$ and an oriented separation $(A_Z,B_Z)$ of~$G[Z]$. Then declare as open the set $O(A_Z,B_Z)$ of all $(A,B)\in\vS$ such that $A\cap Z = A_Z$ and $B\cap Z = B_Z$. We shall say that these $(A,B)$ {\it induce\/} $(A_Z,B_Z)$ on~$Z$, writing $(A_Z,B_Z)=: (A,B)\restricts Z$, and that $(A,B)$ and $(A',B')$ {\it agree on~$Z$} if $(A,B)\restricts Z = (A',B')\restricts Z$.

It is easy to see that the sets $O(A_Z,B_Z)$ do indeed form the basis of a topology on~$\vS$. Indeed, $(A,B)\in\vS$ induces $(A_1,B_1)$ on~$Z_1$ and $(A_2,B_2)$ on~$Z_2$ if and only if it induces on $Z = Z_1\cup Z_2$ some separation $(A_Z,B_Z)$ which in turn induces $(A_i,B_i)$ on~$Z_i$ for both~$i$. Hence $O(A_1,B_1)\cap O(A_2,B_2)$ is the union of all these $O(A_Z,B_Z)$.%
   \COMMENT{}

\proclaim Example \xxxRay.
{\rm If $G$ is a single ray $v_0 v_1\dots$ with end $\omega$, say, then $\tau=\tau_\omega$~is not closed in~$\vS$. Indeed, $\tau$ contains $(\es,V)$ by Lemma \xxxSmall, and hence does not contain~$(V,\es)$. But for every\vadjust{\penalty-200} finite $Z\sub V$ the restriction $(Z,\es)$ of $(V,\es)$ to~$Z$ is also induced by $(\{v_0,\dots,v_n\}, \{v_n, v_{n+1},\dots\})\in \tau$ for every $n$ large enough that $Z\sub\{v_0,\dots,v_{n-1}\}$. So $(V,\es)\in\vS\sm\tau$ has no open neighbourhood in~$\vS\sm\tau.\quad \square$}

\proclaim Example \xxxUFnotClosed.
{\rm Ultrafilter tangles $\tau\in\Theta$ are never closed in~$\vS$. Indeed, let $X\in\X$ witness that $\tau$ is an \uf\ tangle, pick $\C\in U(\tau, X)$, and consider $(A,B)\in\vS$ for $A = X\cup\bigcup\C$ and $B = V\sm\bigcup\C$. This is a separation in~$\vS\sm\tau$ (cf.\ Lemma~\xxxDirect), but every open neighbourhood of $(A,B)$ meets~$\tau$: for every finite $Z\sub V$ we can find a separation $(A',B')\in\tau$ such that $(A',B')\restricts Z = (A,B)\restricts Z$. Such a separation $(A',B')$ can be obtained from $(A,B)$ by moving all the components of~$\C$ that lie in $A\sm Z$ to~$B$.%
   \COMMENT{}
   Then $(A',B')\in\tau$, again by Lemma~\xxxDirect\ (and Lemma~\xxxSymDiff). The details are left to the reader; Theorem~\xxxBlocks\ below includes a formal proof.}
   \endproof

Are any \AT s closed in~$\vS$? As we have seen, they must be end tangles. And such end tangles do exist. Here is the example promised earlier:

\proclaim Example \xxxComplete.
{\rm If $K\sub V$ spans an infinite complete graph in~$G$, then the \AT{}
 $$\tau = \{\,(A,B)\in\vS\mid K\sub B\,\}\eqno\eqK$$
 is closed in~$\vS$. We omit the easy proof.}%
   \COMMENT{}
   \endproof

Perhaps surprisingly, it is not hard to characterize the \AT s that are closed. They are all essentially like Example~\xxxComplete: we just have to generalize the infinite complete subgraph used appropriately. Of the two obvious generalizations, infinite complete minors~[\the\refRSTminors] or subdivisions of infinite complete graphs~[\the\refRST], the latter turns out to be the right one.

Let $\kappa$ be any cardinal. A~set of at least $\kappa$ vertices of $G$ is {\it $(<\kappa)\,$-\th inseparable\/} if no two of them can be separated in~$G$ by fewer than $\kappa$ vertices. A~maximal ${(<\kappa)}\,$-\th inseparable set of vertices is a {\it $\kappa$-block\/}. For example, the branch vertices of a~$TK_\kappa$ are $(<\kappa)\,$-\th inseparable. Conversely:

\proclaim Lemma~\xxxTK.
When $\kappa$ is infinite, every $(<\kappa)\,$-\th inseparable set of vertices in~$G$ contains the branch vertices of some $TK_\kappa\sub G$.

\proof
Let $K\sub V$ be $(<\kappa)\,$-\th inseparable. Viewing $\kappa$ as an ordinal we can find, inductively for all~$\alpha < \kappa$, distinct vertices $v_\alpha\in K$ and internally disjoint $v_\alpha$--$v_\beta$ paths in $G$ for all~$\beta < \alpha$ that also have no inner vertices among those~$v_\beta$ or on any of the paths chosen earlier; this is because $|K|\ge\kappa$, and no two vertices of~$K$ can be separated in~$G$ by the $<\kappa$ vertices used up to that time.
  \endproof

   \COMMENT{}

\noindent
   Note, however, that a $\kappa$-block in~$G$ need not itself be the set of branch vertices of a $TK_\kappa$, even if it has size exactly~$\kappa$. (Consider a~$K_\kappa$ minus an edge.)

\medbreak

We can now prove our last remaining theorem. Let us say that a set ${K\sub V}$\penalty-200\ {\it defines\/} an \AT~$\tau$ if $\tau$ satisfies~\eqK. Note that all the infinite subsets of an $\aleph_0$-block define the same $\aleph_0$-tangle.%
   \COMMENT{}

\proclaim Theorem~\xxxBlocks.
Let $G$ be any graph.
\pitem{i} $\!$The \AT s in $G$ that are not end tangles are never closed in~$\vS$.
\pitem{ii} $\!$An end tangle in $G$ is closed in~$\vS$ if and only if it is defined by an $\A$-block, or equivalently, by the set of branch vertices of some~$TK_\A$.%
   \COMMENT{}
   \enditem

\proof
For every $TK_\A = H\sub G$ there is a unique end~$\omega$ of~$G$ containing all the rays in~$H$. Then $\tau\in\Theta$ is defined by the set of branch vertices of this~$TK_\A$ if and only if $\tau=\tau_\omega$. In view of Lemma~\xxxTK\ it thus suffices to show that an arbitrary $\tau\in\Theta$%
   \COMMENT{}
   is closed in~$\vS$ if and only if it is defined by an $\A$-block.%
   \COMMENT{}

\penalty-2000

Suppose first that $\tau$ is defined by an $\A$-block~$K$. To show that $\tau$ is closed, we have to find for every $(A,B)\in\vS\sm\tau$ a%
   \COMMENT{}
   finite set $Z\sub V$ such that no $(A',B')\in\vS$ that agrees with $(A,B)$ on~$Z$ lies in~$\tau$. As $(A,B)\notin\tau$, we have $K\sub A$; pick $z\in K\sm B$. Then every $(A',B')\in\vS$ that agrees with $(A,B)$ on $Z:=\{z\}$%
   \COMMENT{}
   also also lies in~$\vS\sm\tau$, since $z\in A'\sm B'$ and this implies $K\!\not\sub B'$.

Conversely, consider any $\tau\in\Theta$%
   \COMMENT{}
   and let
 $$K:= \bigcap\{\,B\mid (A,B)\in\tau\,\}.$$
  No two vertices in~$K$ can be separated by in~$G$ by a finite-order separation: one orientation $(A,B)$ of this separation would be in~$\tau$, which would contradict the definition of~$K$ since $A\sm B$ also meets~$K$. If $K$ is infinite, it will clearly be maximal with this property,%
   \COMMENT{}
   and hence be an $\A$-block.%
   \COMMENT{}
   This $\A$-block $K$ will define~$\tau$: by definition of $K$ we have $K\sub B$ for every $(A,B)\in\tau$, while also every $(A,B)\in\vS$ with $K\sub B$ must be in~$\tau$: otherwise $(B,A)\in\tau$ and hence $K\sub A$ by definition of~$K$, but $K\not\sub A\cap B$ because this is finite. Hence $\tau$ will be defined by an $\A$-block, as desired for the forward implication.%
   \Footnote{Whether or not $\tau$ is closed in~$\vS$ is immaterial; we just did not use this assumption.}

It thus suffices to show that if $K$ is finite then $\tau$ is not closed in~$\vS$, which we shall do next.

Assume that $K$ is finite. We have to find%
   \COMMENT{}
   some $(A,B)\in\vS\sm\tau$ that is a limit point of~$\tau$, i.e., which agrees on every finite $Z\sub V$ with some $(A',B')\in\tau$. We choose $(A,B):= (V,K) \in\vS\sm\tau$ (Lemma \xxxSmall).%
   \COMMENT{}

To complete our proof as outlined, let any finite set $Z\sub V$ be given. For every $z\in Z\sm K$%
   \COMMENT{}
   choose $(A_z,B_z)\in\tau$ with $z\in A_z\sm B_z$: this exists, because ${z\notin K}$.%
   \COMMENT{}
   By Lemma~\xxxP, the supremum of all these elements of~$\tau$ and $(K,V)\in\tau$ is again in~$\tau$: we have $(A',B')\in\tau$ for
 $$A':= K\cup\bigcup_{z\in Z\sm K} A_z\quad{\rm and}\quad B':= V\cap\bigcap_{z\in Z\sm K} B_z\,.$$
 As desired, $(A',B')\restricts Z = (A,B)\restricts Z$ (which is $(Z,Z\cap K)$, since $(A,B) = (V,K)$): every $z\in Z\sm K$ lies in some~$A_z$ and outside that~$B_z$, so $z\in A'\sm B'$, while every $z\in Z\cap K$ lies in $K\sub A'$ and also, by definition of~$K$, in every $B_z$ (and hence in~$B'$), since $(A_z,B_z)\in\tau$.
   \endproof

\beginsection \secOutlook. Outlook

There are some obvious leads the reader may like to follow up, as well as one not so obvious one. 

The most obvious is to study the space $|G|$ more closely. There are plenty of basic questions about $|G|$ that we have not even addressed.%
   \COMMENT{}
   For example, how is $|G|$ related to the \SC\ compactification of~$G$?%
   \COMMENT{}
   For which $G$ is $|G|$ the coarsest compactification in which its ends appear as distinct points?%
   \COMMENT{}
   If it is not,%
   \COMMENT{}
   is there a unique such topology, and is there a canonical way to obtain it from~$|G|$?

More important, and probably a good guidance also for which of these basic questions to address, is the potential of $|G|$ for applications in graph theory. For locally finite graphs,\vadjust{\penalty-200} the study of its end compactification $|G|$ has proved very enlightening indeed, and has led to some considerable advances even for purely graph-theoretic problems not originally involving ends~[\the\refTopSurvey]. Might considering our tangle compactification $|G|$ lead to similar advances for arbitrary infinite graphs~$G$?

Another obvious lead is to consider $\kappa$-tangles for cardinals $\kappa>\A$. Do the $\kappa$-tangles that are closed in the space $\vS_\kappa$ of all oriented separations of~$G$ of order~$<\kappa$ form interesting highly connected substructures that do not coincide with classical such structures such as $TK_\kappa$ subgraphs?

Finally, there is an intriguing way to generalize \AT s to separation systems of much more general discrete structures than graphs, introduced in~[\the\refASS]. Essentially, all we need to remember of~$\vS$ is that it is a poset with an order-reversing involution. One can then define stars of `oriented separations' (elements of~$\vS$) as earlier in Section~\secTangles, and for a set $\F$ of such stars one can consider $\F$-tangles. Perhaps there is a natural (submodular) `order' function on~$\vS$, as is the case, for example, for separations in matroids. But even if not, there is a way of expressing \AT s in this framework without any reference to an order function~-- or, indeed, to the cardinality of $\bigcap_i B_i$ as in the definition of~$\T_{<\aleph_0}$.\looseness=-1

We need one more definition to express this. Call an oriented separation $\vs\in\vS$ {\it small\/} if $\vs\le\sv$, where $\sv$ denotes the image of $\vs$ under the involution.%
   \Footnote{For $\vs=(A,B)$, this would be $\sv=(B,A)$. The small separations in our $\vS$ are those of the form $(A,V)$: they satisfy $(A,V)\le(V,A)$, and are the only separations with this property.\looseness=-1}
   Using this term, we can rephrase the definition of a $\T_{<\aleph_0}$-tangle of~$S$ without mentioning cardinalities:

\proclaim Observation.  The $\T_{<\aleph_0}$-tangles of $S$ are the consistent orientations $\tau$ of~$S$ such that no finite star $\sigma\sub\tau$ has a supremum in~$\vS$ whose inverse ist small.

\proof A $\T_{<\aleph_0}$-tangle cannot contain such a star $\sigma = \{\,(A_i,B_i)\mid i= 1,\dots,n\,\}$: since the supremum of~$\sigma$ is $(\bigcup_i A_i,\bigcap_i B_i)$, the inverse of this can be small only if $\bigcap_i B_i$ is finite, which would place $\sigma$ in~$\T_{<\aleph_0}$.

Conversely, let us show that if $\tau$ is a consistent orientation of~$S$ such that no star $\sigma\sub\tau$ has a supremum in~$\vS$ with a small inverse, then $\tau$ has no subset in~$\T_{<\aleph_0}$. For let $\sigma = \{\,(A_i,B_i)\mid i= 1,\dots,n\,\}\sub\tau$ be such a subset. Then $X:=\bigcap_i B_i$ is finite, and the separations $(A'_i,B_i)\ge (A_i,B_i)$ with $A'_i:= A_i\cup X$ still lie in~$\vS$. In fact, they must also lie in~$\tau$. For if $(B_i,A'_i)\in\tau$ for some~$i$, then $\{(A_i,B_i), (B_i,A'_i)\}\sub\tau$ is a star whose supremum $(A_i\cup B_i, B_i\cap A'_i) = (V,X)$ has a small inverse. But the supremum of $\sigma' := \{\,(A'_i,B_i)\mid i=1,\dots,n\,\}\sub\tau$ is $(\bigcup A'_i, \bigcap B_i) = (V,X)$, which has a small inverse~-- a contradiction to the choice of~$\tau$.
   \endproof

If our characterization of the \AT s of $G$ in terms of~$\U$ can be re-done in this abstract setting, it may become meaningful to consider \uf\ tangles in more general structures than graphs, such as ends in matroids, that have been sought for some time.

\beginsection{Acknowledgement}

Thanks to Jakob Kneip for pointing out an oversight in the published version of Lemma~\xxxTK; compare the remark now following the lemma. The proof of Theorem~\xxxBlocks\ uses only the corrected version now stated and proved in Section~\secBlocks.

\beginsection References

\ref\refPrimenden
   Carath\'eodory, \"Uber die Begrenzung einfach zusammenh\"angender Gebiete, \MA73 (1913), 323--370.

\ref\refCarmesinEF
   J.\th Carmesin, All graphs have tree-decompositions displaying their topological ends, \Comb39 (2019), 545--596.

\ref\refCWoess
   D.I.\th Cartwright, P.M.\th Soardi and W.\th Woess, Martin and end compactifications for non-locally-finite graphs, \TAMS338 (1993), 679--693.

\ref\refTopSurvey
   R.\th Diestel, Locally finite graphs with ends: a topological approach. {\sl Discrete Mathematics, Special Volumes \bf 311--312} (2010--11), arXiv:0912.4213.

\ref\refBook
   R.\th Diestel, {\it Graph Theory\/} (5th edition), Springer-Verlag Heidelberg, 2017.

\ref\refASS
   R.\th Diestel, Abstract separation systems, \Order35 (2018), 157--170.

\ref\refTopEnds
   R.\th Diestel \& D.\th K\"uhn, Graph-theoretical versus topological ends of graphs, \JCTB87 (2003), 197--206.

\ref\refTangleTreeAbstract
   R.\th Diestel \& S.\th Oum, Tangle-tree duality in abstract separation systems, \Advances377 (2021), 107470.

\ref\refDualityGraphsMatroids
   R.\th Diestel \& S.\th Oum, Tangle-tree duality: in graphs, matroids and beyond, \Comb39 (2019), 879--910.

\ref\refFreudenthal
      H.\th Freudenthal, \"Uber die Enden topologischer R\"aume
      und Gruppen, \MZ33 (1931), 692--713.

\ref\refFreudenthalNeu
H.\th Freudenthal, Neuaufbau der Endentheorie, \Annals43 (1942), 261--279.

\ref\refHalinEnds
R.\th Halin, \"Uber unendliche Wege in Graphen, \MA157 (1964), 125--137.

\ref\refHalinGrid
R.{\th}Halin, \"Uber die Maximalzahl frem\-der unendlicher Wege, \MN30 (1965), 63--85.

\ref\refHeuerFullGrid
K.\th Heuer, Excluding a full grid minor, \Abh87 (2017), 265--274.

\ref\refPolat
     N.\th Polat, Topological aspects of infinite graphs, in:
     (G.\th Hahn et al., Eds.) {\sl Cycles and Rays\/}, NATO ASI Ser.\th C,
     Kluwer Academic Publishers, Dordrecht 1990.

\ref\refInvLim
L.\th Ribes and P.\th Zalesskii, {\it Profinite Groups\/}, Springer 2010.

\ref\refGMX
   N.\th Robertson \& P.D.\th Seymour, Graph minors~{X}, {O}bstructions to tree-decomposition, \JCTB52 (1991), 153--190.

\ref\refRST
       N.\th Robertson, P.\th D.\th Seymour and R.\th Thomas,
       Excluding subdivisions of infinite cliques, \TAMS332 (1992), 211--223.

\ref\refRSTminors
       N.\th Robertson, P.\th D.\th Seymour and R.\th Thomas,
       {\it Excluding Infinite Clique Minors\/}, {\sl Memoirs of the AMS~\bf 118}, American Mathematical Society, 1995.

\ref\refCTW
C.\th Thomassen and W.\th Woess, Vertex-transitive graphs and accessibility, \JCTB58 (1993), 248--268.

\bigskip\vfill
\ninepoint\noindent Version 27.2.2021
\eject\end